\documentclass[10pt,a4paper,oneside]{amsart}

\usepackage{amssymb,amscd}
\usepackage[latin1]{inputenc}
\usepackage{graphicx}

\def\N{\mathbb{N}}
\def\Z{\mathbb{Z}}

\def\ZZ{\Z \oplus \Z}

\def\P{\pi_{1}}
\def\a{\alpha}
\def\b{\beta}
\def\g{\gamma}
\def\d{\delta}
\def\e{\varepsilon}
\def\f{\phi}
\def\vf{\varphi}

\def\s{\sigma}
\def\w{\omega}
\def\G{\Gamma}

\def\cal{\mathcal}

\def\el{\'el\'ement}

\newtheorem{thm}{Th\'eor\`eme}[section]
\newtheorem{cor}{Corollaire}[section]

\newtheorem{prop}{Proposition}[section]
\newtheorem{defn}{D\'efinition}[section]

\makeindex \makeglossary
\title{Centre, commutativité et conjugaison dans un graphe de groupe}
\begin{document}
\maketitle
\begin{center}
{\sc Jean-Philippe PR\' EAUX}\footnote[1]{\noindent Centre de
recherche de l'école de l'air, Ecole de l'air, F-13661 Salon de
Provence air,}
\footnote[2]{C.M.I., Universit\'e de Provence, 39 rue
F.Joliot-Curie, F-13453 marseille cedex 13.
\\ \indent\ {\it E-mail }: \ preaux@cmi.univ-mrs.fr\smallskip }
\end{center}

\begin{abstract}
We give  characterizations of the center, of  conjugated and of
commuting elements in a fundamental group of a graph of group. We
deduce various results : on the one hand we give a sufficient
condition for the center, the centralizers, and the root
structures in such a group to be in some sense trivial, and on the
other hand we  prove that for any group $G$, the conjugacy problem
reduces to the same problem in a double of $G$
along any finite family of subgroups.\medskip\\
{\sc Résumé.} Nous caractérisons le centre, les éléments conjugués
et les éléments qui commutent dans le groupe fondamental d'un
graphe de groupe. Nous en déduisons divers résultats : d'une part
nous donnons une condition suffisante pour que le centre, les
centralisateurs et la structure des racines d'un tel groupe soient
dans un sens triviaux, et d'autre part nous montrons que pour un
groupe $G$ quelconque  le problème de conjugaison se réduit au
même problème dans un double de $G$ le long d'une famille finie
quelconque de sous-groupes.
\end{abstract}

\section{Introduction}
Amalgames et extension HNN sont des constructions de théorie des
groupe  amplement étudiées et fécondes en exemples (\cite{mks,
rot, Lyndon}. Serre les a généralisé par le concept de groupe
agissant sans inversion sur un arbre (nous parlerons ici plutôt de
groupe fondamental de graphe de groupe, ou plus abusivement de
graphe de groupe) (\cite{serre, dd}).

Nous généralisons ici des résultats connus  pour les amalgames au
cadre le plus général d'un graphe de groupe. Plus précisément nous
donnons une caractérisation  du centre, des éléments qui commutent
et des éléments conjugués dans le groupe fondamental d'un graphe
de groupe (théorèmes \ref{graph center}, \ref{graph_com}  et
\ref{graph_conj}), et en particulier dans une extension HNN
(théorèmes \ref{hnnconj2}, \ref{hnncom} et \ref{hnncenter}).

Nous mettons ensuite à profit l'étude effectuée pour établir
divers résultats. Nous donnons une condition combinatoire
suffisante pour que dans un graphe de groupe les centralisateurs,
le centre et la structure des racines soient dans un sens triviaux
(théorèmes \ref{centralizer}, \ref{triv center} et \ref{srt}) ;
bien que loin d'être nécessaire cette condition est assez large
pour englober de nombreux exemples. Nous montrons aussi que donné
un groupe $G$ quelconque et une famille finie quelconque de
sous-groupes, le problème de conjugaison dans $G$ se réduit au
problème de conjugaison dans le double de $G$ le long de ces
sous-groupes (théorème \ref{double}).

Ce travail est pour une part significative extrait de ma thèse de
doctorat (\cite{phd}) ; je tiens à remercier chaleureusement mon
directeur de thèse Hamish Short pour sa confiance et son soutien
durant toutes ces années.\smallskip
\\
\centerline{Notations}

Nos preuves procédant souvent par disjonction de cas, nous
utiliserons dans une démonstration le symbôle $\square$ pour
dénoter la fin d'un cas, tandis que le symbôle $\blacksquare$
dénotera la fin de la preuve.

\section{Rappels sur les produits amalgam\'es et extensions HNN}
 Nous rappelons les d\'efinitions et les
propri\'et\'es fondamentales d'un produit amalgam\'e, et d'une
extension HNN. Nous ne ferons qu'\'enoncer les r\'esultats sans en
fournir de preuve, le lecteur pouvant se r\'ef\'erer aux ouvrages
\cite{mks}, \cite{Lyndon}, \cite{rot}.\medskip\\
\noindent\textbf{Produits amalgam\'es.}
Soient $A$ et $B$ des groupes donn\'es par les pr\'esentations
respectives \mbox{$< S\;| \; R>$} et $<S'\; | \; R'>$. Soient
$C_A$ un sous-groupe de $A$, $C_B$ un sous-groupe de $B$, et $\f :
C_A \longrightarrow C_B$ un isomorphisme. On appelle
\textsl{produit amalgam\'e} de $A$ et $B$ le long de $\f$, le
groupe que l'on note $A\,\ast_\f B$, donn\'e par la pr\'esentation
:
$$A\,\ast_\f B \ \cong\ < S \cup S' \; | \; R \cup R' \cup
\{\f (c) = c \quad \forall c \in C_A\}>$$ \label{A*phiB}

Si les pr\'esentations de $A$ et $B$ sont finies, et si $C_A$ est
de type fini, $A\,\ast_\f B$ est de pr\'esentation finie. En fait
si $c_1,\ldots ,c_n$ est une famille g\'en\'eratrice de $C_A$, on
a la pr\'esentation finie :
$$ A\,\ast_\f B\ \cong\
< S \cup S' \; | R \cup R' \cup \{\f (c_i) = c_i \quad i=1,\ldots
,n\}>$$

On consid\`ere l'application d'inclusion de $S$ dans $S\cup S'$.
Elle s'\'etend naturellement en un unique homomorphisme du groupe
libre $F(S)$, dans le groupe libre $F(S\cup S')$, qui passe au
quotient pour donner un homomorphisme de $A$ dans $A\,\ast_\f B$.
De la m\^eme fa\c con on d\'efinit un homomorphisme naturel de $B$
dans $A\,\ast_\f B$, qui \'etend l'inclusion de $S'$ dans $S\cup
S'$.

\begin{prop}\label{amalg1}
Les homomorphismes naturels de $A$ et de $B$ dans $A\,\ast_\f B$,
sont injectifs.
\end{prop}

\noindent On les appelera les plongements naturels, et on les
notera respectivement $i_A : A\longrightarrow A\,\ast_\f B$ et
$i_B : B\longrightarrow A\,\ast_\f B$. L'image de $A$ est
appel\'ee premier facteur, celle de $B$ deuxi\`eme facteur. On ne
distinguera pas en g\'en\'eral $A$ et $B$ de leur image $i_A(A)$
ou $i_B(B)$. Ainsi $A$ et $B$ seront consid\'er\'es comme des
sous-groupes de $A\,\ast_\f B$. De plus, les sous-groupes
$i_A(C_A)$ et $i_B(C_B)$ de $A\,\ast_\f B$, sont confondus, en un
sous-groupe que l'on note $C$ ; en fait :
$C=i_A(C_A)=i_B(C_B)=i_A(A)\cap i_B(B)$. Nous pourrons dire, de
fa\c con quelque peu impr\'ecise, que $A\, *_\f B$ est un amalgame
de $A$ et $B$, le long du sous-groupe $C$, et le noter $A\,
*_{\scriptscriptstyle C} B$. \label{A*cB}
\medskip

Si $g$ est un \el\ de  $A\,\ast_\f B$, une \textsl{forme normale}
 pour $g$, est une
suite finie $(g_1,g_2,\ldots ,g_n)$ d'\el s de $A\,\ast_\f B$,
v\'erifiant les conditions :\medskip\\
-- $\forall\, i=1,\ldots ,n$, $g_i\in A \setminus C$, ou
$g_i\in B \setminus C$.\\
-- $\forall\, i=1,\ldots ,n-1$, si $g_i\in A$, alors $g_{i+1}\in B$.\\
-- $g=g_1g_2\cdots g_n$ dans $A\,\ast_\f B$.\medskip\\
\noindent Puisque $A\,\ast_\f B$ a pour famille g\'en\'eratrice
$S\cup S'$, o\`u $S$ engendre $A$, $S'$ engendre $B$ et que $A\cap
B=C$ dans $A\, \ast_\f B$, alors clairement, tout \'el\'ement
admet une forme normale. Dans un certain sens, cette forme normale
est unique. C'est le th\'eor\`eme fondamental des produits libres
amalgam\'es.

\begin{thm}[Forme normale]
\label{amalg2} Consid\'erons un \'el\'ement $g$ de $A\,\ast_\f B$,
ainsi que $(g_1,g_2,\ldots ,g_n)$ et $(g_1',g_2',\ldots ,g_m')$
deux formes normales de $g$. Alors $m=n$,  $\forall\, i=1,\ldots
,n$, $g_i$ et $g_i'$ sont dans un m\^eme facteur, et $g_i^{-1}g_i'
\in C$.
\end{thm}

\label{long1} L'entier $n$ sera appel\'e \textsl{longueur} de
l'\el\ $g$, et not\'e $|g|$. En particulier $|g|=1$ si et
seulement si $g$ est dans un des facteurs.

Nous utiliserons plut\^ot la formulation suivante :
 si $(g_1,g_2,\ldots ,g_n)$ est une forme normale pour $g$,
o\`u les $g_i$ sont des mots sur $S$ ou sur $S'$, alors le mot
$g_1g_2\cdots g_n$ repr\'esente $g$ dans $A\,\ast_\f B$, et pourra
\^etre d\'enomm\'e d'\textsl{\'ecriture sous forme
r\'eduite}.\medskip

\noindent\textbf{Extension HNN.} Soit $A$ un groupe, $C_{-1}$ et
$C_{+1}$ des sous-groupes de $A$, et $\f : C_{-1}\longrightarrow
C_{+1}$ un isomorphisme. Si A admet pour pr\'esentation $<S\ |\
R>$, on appelle \textsl{extension HNN} de A, relativement \`a
$\f$, le groupe que l'on notera $A\ast_\f$, de pr\'esentation,
$$A\ast_\f\ \cong\
<S\cup\{t\}\ |\ R\cup\{t^{-1}ct=\f (c)\quad \forall c\in
C_{-1}\}>$$ \label{A*phi}
 Si le sous-groupe $C_{-1}$ admet une
famille g\'en\'eratrice finie, $c_1,\ldots ,c_n$, alors $A*_\f$
admet la pr\'esentation finie,
$$
A\ast_\f\ \cong\ < S\cup \{t\}\;\mid\;R\cup \{t^{-1}c_it =
\f(c_i)\quad \forall\, i=1,\ldots ,n\}>$$

On consid\`ere l'application identit\'e de $S$ dans $S$. Elle
s'\'etend de fa\c con unique en un homomorphisme de $F(S)$ dans
$F(S\cup\{t\})$, qui passe au quotient pour donner un
homomorphisme naturel de $A$ dans $A\ast_{\f}$.

\begin{prop}
\label{hnn1} L'homorphisme naturel de $A$ dans $A\ast_{\f}$ est un
plongement.
\end{prop}

On note $i_A: A\longrightarrow A\ast_{\f}$ le plongement naturel.
Nous confondrons $A$ et son image $i_A(A)$ par $i_A$, ainsi nous
verrons $A$ comme un sous-groupe de $A\ast_\f$.

Etant donn\'e un \'el\'ement $\g \in A\ast_{\f}$, une
\textsl{forme normale}  pour $g$, est une suite finie d'\el s de
$A*_\f$, de la forme $(g_0,t^{\e_1},g_1,t^{\e_2},\ldots
,t^{\e_n},g_{n})$,
v\'erifiant :\medskip\\
-- $\forall\, i=0,\ldots ,n$, $g_i\in A$, et
$\forall\, i=1,\ldots , n$, $\e_i=\pm 1$.\\
-- $\forall\, i=1,\ldots ,n-1$, si $\e_{i+1}=-\e_i$,
alors $g_{i}\not\in C_{\e_i}$.\\
-- $g=g_0t^{\e_1}g_1t^{\e_2}\cdots t^{\e_n}g_{n}$
dans $A*_\f$.\medskip\\
Une \textsl{forme r\'eduite}  pour $g$ est un mot sur
$S\cup\{t\}$,
 de la forme
$g_0t^{\e_1}g_1t^{\e_2}\cdots t^{\e_n}g_{n}$, o\`u les $g_i$ sont
des mots sur $S$, et $(g_0,t^{\e_1},g_1,t^{\e_2},\ldots
,t^{\e_n},g_{n})$ est une forme normale pour $g$.

Un mot de la forme $t^\e ct^{-\e}$, o\`u $\e=\pm 1$, et $c\in
C_\e$, est appel\'e un \textsl{pinch}. Si un mot contient un tel
pinch, l'op\'eration consistant \`a remplacer ce sous-mot par le
mot $\f^{-\e} (c)$, ne change pas l'\el\ de $A*_\f$
repr\'esent\'e. Ainsi, clairement tout \el\ admet une forme
normale.

Le th\'eor\`eme suivant, fondamental pour les extensions HNN,
garantit l'unicit\'e d'une forme normale. Il est connu sous le nom
de lemme de Britton.

\begin{thm}[Lemme de Britton]
\label{hnn2}
 L'\'el\'ement $1\in A\ast_{\f}$,
admet pour unique forme normale, la suite (1) .
\end{thm}

On  en déduit facilement que si $g$ admet deux formes normales
$(g_0,t^{\e_1},g_1,t^{\e_2},\ldots ,t^{\e_n},g_{n})$ et
$(g'_0,t^{\mu_1},g'_1,t^{\mu_2},\ldots ,t^{\mu_m},g'_{m})$, alors
n\'ecessairement, $n=m$, $\forall i=1,\ldots ,n$, $\e_i=\mu_i$, et
${g'_0}^{-1}g_0\in C_{-\e_1}$ et $g_{n}{g'_{n}}^{-1}\in C_{\e_n}$.
\label{long2} On appelle \textsl{longueur} de $g$, not\'ee $|g|$,
l'entier d\'efini par $|g|=n+1$. Ainsi, $|g|=1$ si et seulement si
$g\in A$.

\section{Rappels sur les graphes de groupe}
Un \textsl{graphe fini} $X$  est la donn\'ee de deux ensembles
finis $\cal{A}_X$ et $\cal{S}_X$, d'une involution sans point fixe
\mbox{$\mathrm{j} : \cal{A}_X\longrightarrow \cal{A}_X$}, et de
deux applications $\mathrm{o},\mathrm{e}: \cal{A}_X\longrightarrow
\cal{S}_X$, v\'erifiant $\forall a\in \cal{A}_X,
\mathrm{e}(a)=\mathrm{o}\circ\mathrm{j}(a)$. \label{or,ext}

\label{ens,som,ar} Les \el s de $\cal{A}_X$ seront appel\'es des
\textsl{ar\^etes},  et les \el s de $\cal{S}_X$, des
\textsl{sommets}.  On notera $-a=\mathrm{j}(a)$, et les ar\^etes
$a$ et $-a$ seront dites oppos\'ees. Si $a$ est une arête,
$\{a,-a\}$ sera appelé une {\sl arête non-orientée}.

Une \textsl{orientation} de $X$, est un sous-ensemble
$\cal{A}^+_X$ de $\cal{A}_X$, contenant, pour tout
$a\in\cal{A}_X$,
 exactement un \el\
de l'ensemble $\{ a,-a\}$. Un graphe muni d'une orientation sera
dit \textsl{orient\'e}.

Un \textsl{chemin}  de $X$ est une suite finie d'ar\^etes,
$c=(a_1,\ldots ,a_n)$, qui v\'erifie pour tout \linebreak
\mbox{$i=1,\ldots ,p-1$}, $\mathrm{e}(a_i)=\mathrm{o}(a_{i+1})$.
Les sommets $\mathrm{o}(a_1)$ et $\mathrm{e}(a_n)$, sont appel\'es
respectivement origine et extr\'emit\'e du chemin $c$. Le chemin
est dit \textsl{ferm\'e} si $\mathrm{o}(a_1)=\mathrm{e}(a_n)$. Il
est dit \textsl{r\'eduit}, si $\forall i=1,\ldots ,p-1$, $a_i\not=
-a_{i+1}$.

Un graphe $X$ est dit \textsl{connexe},  si pour tout couple de
sommets, $s_1,s_2$, il existe un chemin de $X$ ayant pour origine
$s_1$, et pour extr\'emit\'e $s_2$.

Un graphe est un \textsl{arbre}, si pour tout couple de sommets
$s_1,s_2$, il existe un unique chemin r\'eduit d'origine $s_1$ et
d'extr\'emit\'e $s_2$. C'est un fait bien connu, que tout graphe
fini $X$ contient un arbre maximal, \emph{i.e.} un sous-graphe qui
est un arbre, et qui
a m\^eme sommets que $X$. Il n'est en g\'en\'eral pas unique.\medskip\\
\indent Un \textsl{graphe de groupe} $(\cal{G},X)$ \label{(G,X)}
 est la donn\'ee d'un graphe
fini connexe, orient\'e, et d'une famille $\cal{G}$, consistant en
: un ensemble $\cal{GA}_X=\{G_a, a\in \cal{A}^+_X\}$ de groupes,
appel\'es \textsl{groupes d'ar\^ete}, \label{ens,gpe,ar,som}
un ensemble $\cal{GS}_X=\{G_s, s\in \cal{S}_X\}$ de groupes,
appel\'es \textsl{groupes de sommet}, et pour tout $a \in
\cal{A}^+_X$, les monomorphismes $\f_a :G_a\longrightarrow
G_{\mathrm{o}(a)}$ et $\f_{-a} :G_a\longrightarrow
G_{\mathrm{e}(a)}$. \label{phi_a}

Un \textsl{lacet}  dans $(\cal{G},X)$ est une suite de la forme
$(g_0,a_1,g_1,\ldots ,a_n,g_n)$, o\`u \linebreak $(a_1,a_2,\ldots
,a_n)$ est un chemin ferm\'e de $X$, et $\forall k=1,\ldots ,n\,$,
$g_{k-1}\in G_{\mathrm{o}(a_i)}$, et $g_n\in G_{\mathrm{e}(a_n)}$.
On dit qu'il a pour origine $\mathrm{o}(a_1)$.

Il existe une relation d'\'equivalence $\equiv$ pour les lacets de
$(\cal{G},X)$, engendr\'ee par les relations :
\begin{gather*}
\forall a \in \cal{A}_X, \forall h\in G_a, \forall g\in
G_{\mathrm{o}(a)}\qquad (g,a,\f_{-a}(h), -a, \f_a(h)^{-1})=(g)\\
\forall a\in \cal{A}_X, \forall g,g' \in G_{\mathrm{o}(a)}, \qquad
(g,a,1,-a,g')= (gg')
\end{gather*}
 Si $s_0$ est un sommet de $X$, l'ensemble des
lacets d'origine $s_0$, admet une op\'eration, appel\'ee
op\'eration de concat\'enation :
$$(g_0,a_1,\ldots ,a_p,g_p).(g_{p+1}, a_{p+1},\ldots a_n, g_{n+1})
=(g_0,a_1,\ldots ,a_p,g_pg_{p+1}, a_{p+1},\ldots a_n, g_{n+1})$$
Cette op\'eration passe au quotient sous la relation $\equiv$ en
une op\'eration sur l'ensemble des classes d'\'equivalence de
lacets bas\'es en $s_0$, qui admet d\`es-lors une structure de
groupe. On l'appelle le \textsl{groupe fondamental} de
$(\cal{G},X)$, bas\'e en $s_0$,  et on le note $\P
(\cal{G},X,s_0)$. Si tous les groupes de sommets (et donc les
groupes d'ar\^etes) sont triviaux, on retrouve le groupe
fondamental de $X$ bas\'e en $s_0$.

Le th\'eor\`eme qui suit est un r\'esultat fondamental pour les
graphes de groupe. C'est le seul r\'esultat sur les graphes de
groupe que nous emploierons pour mener \`a bien notre \'etude,
d\'elaissant ici le concept de groupe agissant sans inversion sur
un arbre
 (cf. \cite{serre}, \cite{dd}).

\begin{thm}
\label{graphpres}
Soient $(\cal{G},X)$ un graphe de groupe,  $s_0$ un sommet de $X$,
et $T$ un arbre maximal de $X$.

Alors $\P (\cal{G},X,s_0)$ est le groupe obtenu \`a partir du
produit libre $({\ast}_{s\in \cal{S}_X} G_s)\ast F$, o\`u $F$ est
le groupe libre engendr\'e par $\{t_a\, ;a \in \cal{A}_X\}$, en
ajoutant les relations :
\begin{align*}
\f_a(h)&= t_a\,\f_{-a}(h)\, t_a^{-1} \quad\forall a\in
\cal{A}^+_X,
\forall h\in G_a\\
t_{-a}&=t_a^{-1}\qquad\qquad\quad\forall a\in \cal{A}_X\\
t_a &=1  \qquad\qquad\qquad\forall a\in \cal{A}_X\cap \cal{A}_T
\end{align*}
\label{t_a}
\end{thm}
\noindent \textbf{Remarque 1 :} Ainsi $\P (\cal{G},X,s_0)$ ne
d\'epend pas du choix du point de base $s_0$. On parlera du groupe
fondamental de $(\cal{G},X)$ que l'on notera $\P (\cal{G},X)$.\medskip\\
%
%
\noindent\textbf{Remarque 2 :} Un graphe de groupe $(\cal{G},X)$,
muni d'un arbre maximal $T$ de $X$, est appel\'e un \textsl{graphe
d\'ecompos\'e}.  On pourra le noter $(\cal{G},X,T)$.
\label{(G,X,T)}
 Donn\'e un arbre maximal $T$ de $X$, un \el\
de $\cal{A}^+_T$ sera appel\'e \textsl{ar\^ete $T$-s\'eparante} ;
un \el\ de $\cal A^+_X\cap (\cal A_X\setminus \cal A_T)$ sera
appelé {\sl arête non $T$-séparante}.
\medskip
\\
\noindent \textbf{Remarque 3 :} Si $X=T$ est un arbre, $\P
(\cal{G},X)$ est l'amalgame des groupes de sommets $G_s, s\in
\cal{S}_T$, le long des isomorphismes $\f_{-a}\circ \f_{a}^{-1}$
d\'efinis pour tout $a\in \cal{A}_T^+$. Ainsi pour tout $s\in
\cal{S}_T$, on a le morphisme injectif $p_s : G_s \hookrightarrow
\P(\cal{G},T)$

Si $X$ n'est pas un arbre,  $\P(\cal{G},X)$ est l'extension HNN de
$\P(\cal{G},T)$, le long des isomorphismes
$p_{\mathrm{e}(a)}\circ\f_{-a}\circ \f_{a}^{-1}\circ
p_{\mathrm{o}(a)}^{-1}$ pour toute ar\^ete $a$ non
$T$-s\'eparante. Ainsi on a le morphisme injectif $\Pi
:\P(\cal{G},T)\hookrightarrow \P(\cal{G},X)$

On a donc, pour tout sommet $s\in \cal{S}_T=\cal{S}_X$, un
plongement naturel $\Pi_s=\Pi\circ p_s$ de $G_s$ dans
$\P(\cal{G},X)$. Ces plongements ne d\'ependent que de la donn\'ee
d'un arbre maximal $T$ de $X$. Lorsque $T$ sera fix\'e (ce qui
sera toujours le cas), on commettra l'abus de langage de confondre
$\Pi_s(G_s)$ avec $G_s$. Ainsi on pourra voir $G_s$ comme un
sous-groupe de $\P(\cal{G},X)$. Un tel sous-groupe de
$\P(\cal{G},X)$ sera appel\'e \textsl{sous-groupe de sommet}.

Pour toute ar\^ete $a$ $T$-s\'eparante d'origine $s_1$ et
d'extr\'emit\'e $s_2$, les applications de $G_a$ dans
$\P(\cal{G},X)$, $\Pi_{s_1}\circ \f_a$ et $\Pi_{s_2}\circ \f_{-a}$
sont \'egales. On confondra $G_a$ et $\Pi_{s_1}\circ
\f_a(G_a)=\Pi_{s_2}\circ \f_{-a}(G_a)\subset G_{s_1}\cap
G_{s_2}\subset \P(\cal{G},X)$. \label{G_a}

Soit $a\in \cal{A}_X$ une ar\^ete, d'origine $s_1$ et
d'extr\'emit\'e $s_2$, on notera $G_a^-=\Pi_{s_1}\circ \f_a(G_a)$
et $G_a^+=\Pi_{s_2}\circ \f_{-a}(G_a)$, ainsi que $\varphi_a$,
l'isomorphisme  de $G_a^-$ dans $G_a^+$, d\'efini par
$\varphi_a(h)= \Pi_{s_2}\circ \f_{-a}\circ\f_a^{-1}\circ
\Pi_{s_1}^{-1}(h)$ pour tout $h\in G_a^-$. \label{vf_a}
 Si $a$ est
$T$-s\'eparante, les applications $\Pi_{s_1}\circ \f_a$ et
$\Pi_{s_2}\circ \f_{-a}$ sont identiques, ainsi $G_a^-=G_a^+=G_a$,
et $\varphi_a=\varphi_{-a}$ est l'identit\'e. Si $a$ est non
$T$-s\'eparante, en g\'en\'eral $G_a^-\not= G_a^+$, et dans
$\P(\cal{G},X)$, on
 a la relation
$\forall h \in G_a^-, \quad h=t_a\, \vf_a(h)\, t_a^{-1}$. Les
sous-groupes $G_a^-,G_a^+$ de $\P(\cal{G},X)$ seront appel\'es
\textsl{sous-groupes d'ar\^ete}. Notons que $G_a^-\subset
G_{\mathrm{o}(a)}=G_{s_1}$,
$G_a^+\subset G_{\mathrm{e}(a)}=G_{s_2}$.\medskip\\
\noindent \textbf{Remarque 4 :} Soit $(\cal{G},X)$ un graphe de
groupe, et $Y$ un sous-graphe orient\'e de $X$ (\emph{i.e.} $X$ et
$Y$ sont d'orientations compatibles). On peut construire un
\textsl{sous-graphe de groupe} $(\cal{G}',Y)$ de $(\cal{G},X)$ en
se restreignant aux groupes de sommet  $\cal{G}'\cal{S}_Y=\{ G_s
\in \cal{GS}_X\, ; \, s \in \cal{S}_Y\}$, aux graphes d'arête
$\cal{G}'\cal{A}_Y=\{ G_a \in \cal{GA}_X\, ; \, a \in \cal{A}_Y\}$
et aux monomorphismes de
 $\cal{G}$, $\f_a$ et
$\f_{-a}$ pour tout $a\in \cal{A}_Y^+\subset \cal{A}_X^+$.
Dans la suite, on notera $(\cal{G},Y)$ au lieu de $(\cal{G}',Y)$.
Observons que $\P(\cal{G},Y)$ se plonge naturellement dans
$\P(\cal{G},X)$.

Soit $(\cal G,X,T)$ un graphe de groupe décomposé. Si $a$ est une
arête non $T$-séparante (resp. $T$ séparante) on note $X_1$ (et
$X_2$) la (les) composante(s) connexe(s) du graphe
$X\setminus\{a,-a\}$. Alors $\pi_1(\cal{G},X)$ se décompose en une
extension HNN (en un amalgame) de $\pi_1(\cal{G},X_1)$ (et
$\pi_1(\cal{G} X_2)$) le long de $\vf_a$. On dira que l'on a {\sl
décomposé $(\cal G,X,T)$ le long de l'arête $a$}.
\medskip\\
\noindent \textbf{Remarque 5 :} Ce th\'eor\`eme fournit
implicitement un pr\'esentation pour le groupe fondamental $\P
(\cal{G},X)$, d'un graphe de groupe d\'ecompos\'e $(\cal{G},X,T)$.
Si de plus les groupes de sommet sont finiment pr\'esent\'es, et
les groupes d'ar\^ete de type fini, alors $\P (\cal{G},X)$ est
finiment pr\'esent\'e. Lorsque les présentations finies et les
familles génératrices sont fixées, une telle pr\'esentation de
$(\cal{G},X)$ sera appel\'ee \textsl{pr\'esentation canonique} de
$\P(\cal{G},X)$.  La famille g\'en\'eratrice donn\'ee par cette
pr\'esentation, sera not\'ee $\cal{G}en(X)$. \label{genX}

\section{Exemple : graphe de groupe et décomposition JSJ}

\subsection{Th\'eor\`eme Jaco-Shalen-Johannson}

Soit $W$ une surface compacte, \`a deux faces, proprement
plong\'ee dans une $3$-vari\'et\'e $M$. On notera $\s_W(M)$ la
$3$-vari\'et\'e obtenue en d\'ecomposant $M$ le long de $W$. Plus
pr\'ecis\'ement, si \label{s_WM}  $\cal{T}_1,\cal{T}_2,\ldots
,\cal{T}_n$ sont les composantes connexes de $W$, on peut trouver
des voisinages r\'eguliers des composantes, $V(\cal{T}_1),\ldots
,V(\cal{T}_n)$, deux \`a deux disjoints, et $\s_W(M)$ est d\'efini
par :
$$\s_W(M)=M-\bigcup_{i=1}^n \,\mathrm{int}(V(\cal{T}_i))$$
La 3--vari\'et\'e $\s_W(M)$ est en g\'en\'eral non connexe. Si
$W\not= \varnothing$, toutes les composantes connexes de $\s_W(M)$
sont \`a bord non vide. Nous pouvons dès-lors \'enoncer le
th\'eor\`eme de d\'ecomposition des vari\'et\'es Haken ferm\'ees
(cf. \cite{js} corollaire V.5.1.):

\begin{thm}[{\bf Jaco-Shalen-Johannson}]
 \label{jsj} Soit $M$
une 3--vari\'et\'e Haken $\partial$-irréductible, fer\-m\'ee ou à
bord torique. Il existe une surface compacte $W$ proprement
plong\'ee dans $M$, incompressible, \`a deux faces, unique \`a
isotopie ambiante de $M$ pr\`es, v\'erifiant les propri\'et\'es
suivantes :
\begin{itemize}
\item[(i)] Les composantes connexes de $W$ sont  des
tores.
\item[(ii)] Chaque composante connexe de
$\s_W(M)$ est soit un fibr\'e de Seifert, soit atoro\"\i dale.
\item[(iii)] $W$ est minimale pour l'inclusion, dans
la classe des surfaces v\'erifiant \emph{(i)} et \emph{(ii)}.
\end{itemize}
\end{thm}

\noindent {\bf Remarque :} Avec le th\'eor\`eme d'hyperbolisation
de Thurston (cf. \cite{thurston}), les composantes connexes de
$\s_W(M)$ admettent soit une fibration de Seifert, soit une
structure hyperbolique de volume fini.

\subsection{D\'ecomposition d'une 3--vari\'et\'e en pi\`eces \'el\'ementaires}

Dans tout ce paragraphe, on note $M$ une 3--vari\'et\'e vérifiant
les hypothèses du théorème \ref{jsj}.

On appelle \textsl{d\'ecomposition}
 de $M$, une surface $W$
incompressible, \`a deux faces, proprement\linebreak plong\'ee
dans $M$, v\'erifiant les conditions (i) et (ii), du th\'eor\`eme
\ref{jsj}. On appelle \textsl{d\'ecomposition JSJ} une
d\'ecomposition de $M$ minimale au sens du théorème \ref{jsj}. On
dira que la d\'ecomposition est \textsl{triviale} si
$W=\varnothing$.

Consid\'erons une d\'ecomposition $W$ de $M$. Il existe une
application canonique, \label{app,idd}  \textsl{l'application
d'identification}, $r : \s_W(M)\longrightarrow M$, qui est telle
que $r^{-1}(W)$ consiste en deux copies hom\'eo\-morphes de $W$
dans $\partial M$, et que $r$ restreinte \`a $r^{-1}(M- W)$ soit
un hom\'eomorphisme.

Notons $\cal{T}_1,\ldots,\cal{T}_i,\ldots, \cal{T}_k$
 les composantes connexes
de $W$. Puisque  $\cal{T}_i$ est \`a deux faces,
$r^{-1}(\cal{T}_i)$ consiste en deux copies hom\'eomorphes de
$\cal{T}_i$ dans $\partial \s_W(M)$ ; on les note arbitrairement
$\cal{T}_i^-$ et $\cal{T}_i^+$. Les restrictions de $r$ \`a
$\cal{T}_i^-$ et $\cal{T}_i^+$ sont des hom\'eomorphismes sur
$\cal{T}_i$. On d\'efinit alors l'hom\'eomorphisme
$f_i:\cal{T}_i^-\longrightarrow \cal{T}_i^+$, qui fait commuter le
diagramme suivant :
\begin{equation*}
\begin{CD}
\cal{T}_i^- @>f_i>> \cal{T}_i^+ \\
@Vr_{| \cal{T}_i^-}VV           @VVr_{| \cal{T}_i^-}V\\
\cal{T}_i @>\operatorname{Id}>> \cal{T}_i
\end{CD}
\end{equation*}
On appelera $f_i$, l'hom\'eomor\-phisme associ\'e \`a $\cal{T}_i$,
et de fa\c con moins pr\'ecise, un tel $f_i$, un
\textsl{hom\'eomor\-phisme associ\'e \`a la d\'ecomposition}.

On note $M_1,\ldots,M_p$, les composantes connexes de $\s_W(M)$ ;
on les appelle les \textsl{pi\`eces \'el\'ementaires}. Ainsi,
$$M\cong (\bigcup_{i=1\ldots p} M_i)_{\diagup f_1,\ldots,f_p}$$

\subsection{Graphe de groupe provenant de la décomposition JSJ}

Soit $M$ une 3--vari\'et\'e Haken, \`a bord torique non vide.
Notons $\cal{T}_1,\ldots ,\cal{T}_q$ les composantes connexes de
$\partial M$. On appelle {\bf syst\`eme p\'eriph\'eral}
 de $\P(M)$, la donn\'ee pour tout
$i=1,\ldots q$, d'un plongement $p_i:\ZZ\hookrightarrow \P(M)$
induit par l'inclusion de $\cal{T}_i$ dans $M$ (il est
d\'etermin\'e par la donn\'ee d'une classe d'homotopie de chemin
du point de base de $\P(M)$ au point de base de $\P(\cal{T}_i)$).

On se donne une vari\'et\'e Haken $M$ v\'erifiant $\chi(M)=0$, et
une d\'ecomposition $W$ de $M$. On note $\cal{T}_1,\ldots
\cal{T}_q$ les composantes de $W$, $M_1,\ldots M_p$ les pièces
\'el\'ementaires, et $f_1,\ldots ,f_q$ les hom\'eomorphismes
associ\'es \`a la d\'ecomposition. On consid\`ere aussi un
syst\`eme p\'eriph\'eral de $\P(M)$. On note pour $i=1,\ldots ,q$,
le plongement de $\P(\cal{T}_i^-)$ dans $\P(M)$,
$p_i^-:\ZZ\hookrightarrow \P(M)$ et le \linebreak plongement de
$\P(\cal{T}_i^+)$ dans $\P(M)$, $p_i^+:\ZZ\hookrightarrow \P(M)$
et l'on note ${f_i}_\ast$ l'application de
$p_i^-(\ZZ)$ sur $p_i^+(\ZZ)$ induite par $f_i$.\\
On construit un \textsl{graphe de groupe associ\'e}  \`a la
d\'ecomposition $W$ de $M$, de la fa\c con suivante :\\
Le graphe orient\'e $X$ est d\'efini par : $\cal{S}_X=\{
S_1,\ldots ,S_j,\ldots S_p\}$ et $\cal{A}^+_X=\{a_1,\ldots
,a_i,\ldots a_q\}$ avec l'origine de $a_i$, $\mathrm{o}(a_i)=S_r$
o\`u $\cal{T}_i^-\subset M_r$, et l'extr\'emit\'e de $a_i$,
$\mathrm{e}(a_i)=S_t$ o\`u $\cal{T}_i^+\subset M_t$.\\
Les groupes de sommets et les groupes d'ar\^ete sont d\'efinis
respectivement par : $G_{S_j}=\P(M_j)$, $j=1,\ldots ,p$ et
$G_{a_i}=\ZZ$, $i=1,\ldots ,q$ et les monomorphismes par
$\f_{a_i}=p_i^-$, et $\f_{-a_i}={f_{a_i}}_\ast\circ p_i^-$. Alors,
avec les notations de la remarque 3 du th\'eor\`eme
\ref{graphpres}, $G_{a_i}^-=p_i^-(\ZZ), G_{a_i}^+=p_i^+(\ZZ)$, et
$\vf_{a_i}={f_{a_i}}_\ast$. On a bien construit un graphe de
groupe, $(\cal{G},X)$. Ce graphe d\'epend de $W$, et du choix d'un
syst\`eme p\'eriph\'eral. La classe d'isomorphisme de son groupe
fondamental n'en d\'epend pas.

Le r\'esultat suivant, s'obtient par une suite d'applications du
 th\'eor\`eme de Van-Kampen.
\begin{thm}
Soit $M$ une vari\'et\'e Haken, et $W$ une d\'ecomposition. Le
groupe fondamental du graphe de groupe associ\'e $(\cal{G},X)$,
est isomorphe au groupe fondamental de $M$.
\end{thm}

\begin{figure}[ht]
\centerline{\includegraphics[scale=0.5]{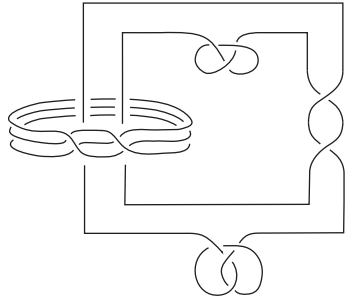}}
\caption{L'entrelacs $L$ de $S^3$} \label{fig1}
\end{figure}

\begin{figure}[!ht]
\centerline{\includegraphics[scale=0.5]{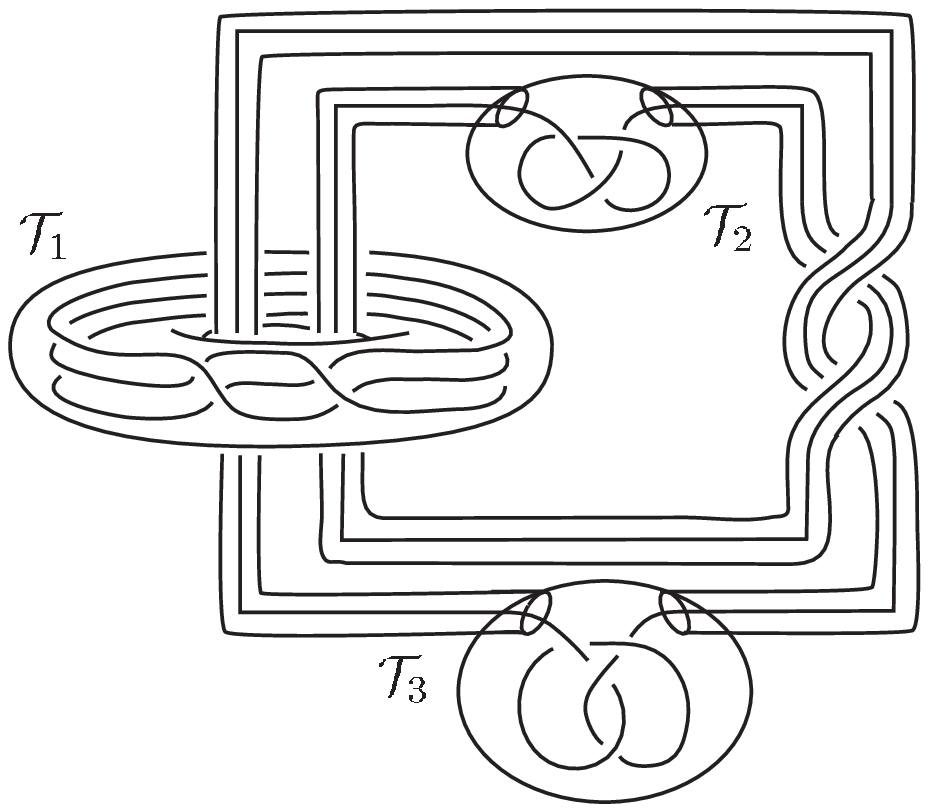}} \caption{Les
tores essentiels $\cal{T}_1$, $\cal{T}_2$, et $\cal{T}_3$ de
$S^3-\mathrm{int}(N(L))$}\label{fig2}
\end{figure}

\begin{figure}[!ht]
\centerline{\includegraphics[scale=0.5]{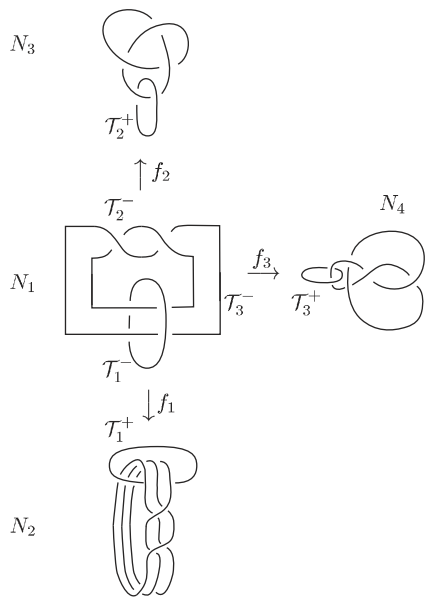}} \caption{Les
pi\`eces \'el\'ementaires $N_1,N_2,N_3,N_4$, et les
hom\'eomorphismes $f_1,f_2,f_3$ associ\'es \`a la
d\'ecomposition}\label{fig3}
\end{figure}

\begin{figure}[ht]
\centerline{\includegraphics[scale=0.5]{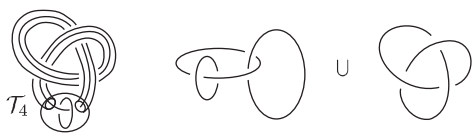}} \caption{Le
tore essentiel $\cal{T}_4$ de $M_2$ (\`a gauche), et les deux
pi\`eces obtenues en d\'ecomposant $M_2$ le long de $\cal{T}_4$
(\`a droite)}\label{fig4}
\end{figure}

\begin{figure}[ht]
\centerline{\includegraphics[scale=0.5]{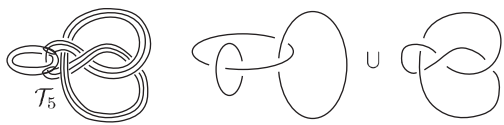}} \caption{Le
tore essentiel de $M_3$ (\`a gauche) et les deux pièces obtenues
en d\'ecomposant le long de $\cal{T}_5$ }\label{fig5}
\end{figure}

\begin{figure}[!ht]
\centerline{\includegraphics[scale=0.5]{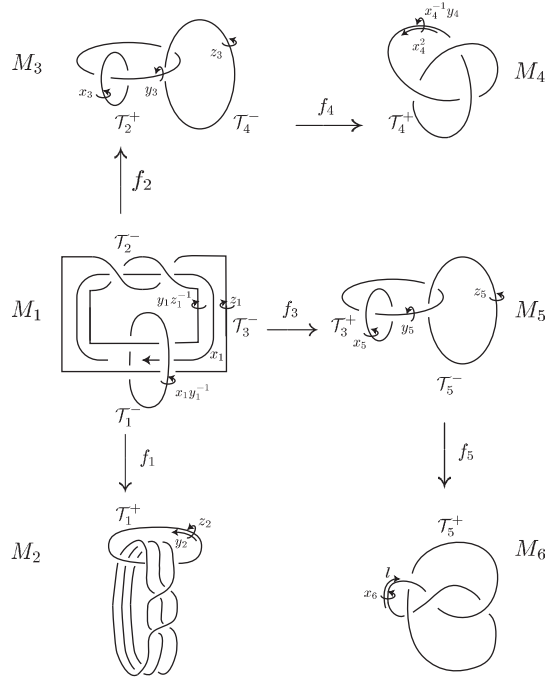}}
\caption{D\'ecomposition de $S^3-\mathrm{int}(N(L))$ le long de
$W$}\label{fig6}
\end{figure}

\noindent\textbf{Exemple :} Consid\'erons l'entrelacs \`a 3
composantes, $L$ de $S^3$ (figure \ref{fig1}).
Consid\'erons trois tores essentiels $\cal{T}_1$, $\cal{T}_2$, et
$\cal{T}_3$ dans le compl\'ement de $L$, $S^3-\mathrm{int}(N(L))$,
donn\'es figure \ref{fig2}.
En d\'ecomposant $S^3-\mathrm{int}(N)(L)$ le long de $\cal{T}_1
\cup \cal{T}_2 \cup \cal{T}_3$, on obtient les vari\'et\'es
$N_1,N_2,N_3,N_4$ (ce sont des compl\'ements d'entrelacs dans
$S^3$), ainsi que des hom\'eomor\-phismes $f_1,f_2,f_3$ associ\'es
\`a la d\'ecomposition (figure \ref{fig3}).
Il est facile de voir que $N_1$ est un fibr\'e en cercle sur la
sph\`ere \`a 3 trous ;  $N_2$ est un fibr\'e de Seifert ayant pour
base la sph\`ere \`a 2 trous, et une fibre exceptionnelle d'indice
$3$ ; et  $N_3$ et $N_4$ contiennent tous deux un tore essentiel.
On note $\cal{T}_4$ le tore essentiel  de $N_3$ (figure
\ref{fig4}), et on d\'ecompose $N_3$ le long de $\cal{T}_4$.
Notons aussi $\cal{T}_5$ le tore essentiel de $N_4$ (figure
\ref{fig5}), et d\'ecomposons $N_4$ le long de $\cal{T}_5$.
Notons $M_1=N_1$, $M_2=N_2$, $M_3$, $M_4$ les pi\`eces obtenues en
d\'ecomposant $N_3$ le long de $\cal{T}_4$ (o\`u $M_4$ est le
tr\`efle), et $M_5,M_6$ les pi\`eces obtenues en d\'ecomposant
$N_4$ le long de $\cal{T}_5$ (o\`u $M_6$ est le noeud de huit).
Les pi\`eces $M_3$ et $M_5$ sont hom\'eomorphes au $S^1$-fibr\'e
trivial  ayant pour base la sph\`ere \`a trois trous. De plus le
compl\'ement du tr\`efle est un fibr\'e de Seifert, et le noeud de
huit est connu pour avoir
 un compl\'ement hyperbolique (de volume fini).
Ainsi, la surface
$W=\cal{T}_1\cup\cal{T}_2\cup\cal{T}_3\cup\cal{T}_4 \cup\cal{T}_5$
est une d\'ecomposition du compl\'ement de l'entrelacs $L$. En
utilisant la classification des compl\'ements d'entrelacs
admettant une fibration de Seifert (\cite{link}), on v\'erifie
imm\'ediatement que $W$ est minimale. Il est facile de se donner
des hom\'eomorphismes $f_1,f_2,f_3$, $f_4$, et $f_5$ associ\'es
\`a la d\'ecomposition (figure \ref{fig6}).
On peut remarquer, bien que cela soit anecdotique, que  $M_1$ et
$M_3$ sont hom\'eomor\-phes, bien qu'ils soient des compl\'ements
d'entrelacs  de types distincts (comparer avec le th\'eor\`eme de
Gordon-Luecke dans le cas d'un noeud).

\begin{figure}[ht]
\centerline{\includegraphics[scale=0.5]{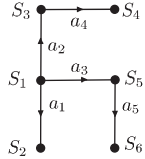}} \caption{Le
graphe orient\'e $X$}\label{fig7}
\end{figure}

 On construit alors le graphe de groupe $(\cal{G},X)$, en consid\'erant tout
d'abord le graphe orient\'e $X$ (figure \ref{fig7}).

puis les groupes de sommets :
$$G_{S_1}=\P(M_1)\cong\ <x_1,y_1,z_1\ |\ [x_1,y_1]=[x_1,z_1]=1>\
\cong F_2\times \Z$$
$$G_{S_2}=\P(M_2)\cong\ <x_2,y_2,z_2\ |\ x_2^3=y_2^2,\ [y_2,z_2]=1>$$
$$G_{S_3}=\P(M_3)\cong\ <x_3,y_3,z_3\ |\ [x_3,y_3]=[y_3,z_3]=1 >
\ \cong F_2\times \Z$$
$$G_{S_4}=\P(M_4)\cong\ < x_4,y_4\ |\ x_4^2=y_4^3>$$
$$G_{S_5}=\P(M_5)\cong\ <x_5,y_5,z_5\ |\ [x_5,y_5]=[y_5,z_5]=1 >
\ \cong F_2\times \Z$$
$$G_{S_6}=\P(M_6)\cong\ <x_6,
y_6\ |\ x_6y_6^{2}x_6^2y_6x_6y_6^{-1}x_6^{-2}y_6^{-1}=1>$$
\normalsize On choisit alors un plongement des groupes des
composantes $\cal{T}_i^{\pm}$ afin d'obtenir les plongements des
groupes d'ar\^etes dans les groupes de sommet ;  les isomorphismes
$\vf_{a_j}$ pour toute ar\^ete $a_j$ sont alors d\'etermin\'es par
les classes d'isotopie des hom\'eomorphismes $f_j$, pour $j$
variant de 1 jusqu'\`a 5. Le calcul fournit :

$$G_{a_1}^-=\P(\cal{T}_1^-)=\ <x_1y_1^{-1},y_1>_{G_{S_1}}\ \subset G_{S_1}$$
$$G_{a_1}^+=\P(\cal{T}_1^+)=\ <z_2,y_2>_{G_{S_2}}\ \subset G_{S_2}$$
$$\vf_{a_1}(x_1y_1^{-1})=y_2\qquad
\vf_{a_1}(y_1)=z_2$$

$$G_{a_2}^-=\P(\cal{T}_2^-)=\ <y_1z_1^{-1},x_1>_{G_{S_1}}\ \subset G_{S_1}$$
$$G_{a_2}^+=\P(\cal{T}_2^+)=\ <x_3,y_3>_{G_{S_3}}\ \subset G_{S_3}$$
$$\vf_{a_2}(y_1z_1^{-1})=y_3\qquad
\vf_{a_2}(x_1)=x_3$$

$$G_{a_3}^-=\P(\cal{T}_3^-)=\ <z_1,x_1>_{G_{S_1}}\ \subset G_{S_1}$$
$$G_{a_3}^+=\P(\cal{T}_3^+)=\ <x_5,y_5>_{G_{S_5}}\ \subset G_{S_5}$$
$$\vf_{a_3}(z_1)=y_5\qquad \vf_{a_3}(x_1)=x_5$$

$$G_{a_4}^-=\P(\cal{T}_4^-)=\ <z_3,y_3>_{G_{S_3}}\ \subset G_{S_3}$$
$$G_{a_4}^+=\P(\cal{T}_4^+)=\ <x_4^{-1}y_4,x_4^2>_{G_{S_4}}\ \subset G_{S_4}$$
$$\vf_{a_4}(z_3)=x_4^2\qquad \vf_{a_4}(y_3)=x_4^{-1}y_4$$

$$G_{a_5}^-=\P(\cal{T}_5^-)=\ <z_5,y_5>_{G_{S_5}}\ \subset G_{S_5}$$
$$G_{a_5}^+=\P(\cal{T}_5^+)=\ <x_6,
l=x_6y_6x_6y_6^{-1}x_6^{-1}y_6^2x_6^2y_6x_6^{-1}>_{G_{S_6}}\
\subset G_{S_6}$$
$$\vf_{a_5}(z_5)=l\qquad \vf_{a_5}(y_5)=x_6$$
Ceci nous suffit \`a d\'efinir le graphe de groupe $(\cal{G},X)$.
Il d\'epend clairement d'un choix des plongements des composantes
au bord. Son groupe fondamental, cependant, n'en d\'epend pas ;
c'est aussi le groupe fondamental du compl\'ement de $L$ dans
$S^3$. Le th\'eor\`eme \ref{graphpres} fournit la pr\'esentation :
\begin{align*}
<x_1,y_1,z_1,x_2,y_2,z_2,x_3,y_3,z_3,x_4,y_4,x_5,y_5,z_5,x_6,y_6\ |&\\
[x_1,y_1]=[x_1,z_1]=[y_2,z_2]=[x_3,y_3]=[y_3&,z_3]=[x_5,y_5]=[y_5,z_5]=1,\\
x_2^3=y_2^2,\quad x_4^2=y_4^3,&\quad
x_6y_6^{2}x_6^2y_6x_6y_6^{-1}x_6^{-2}y_6^{-1}=1,\\
x_1y_1^{-1}=y_2,\quad y_1=&z_2,\quad y_1z_1^{-1}=y_3,\quad
x_1=x_3,\\
 z_1=y_5,\quad
 x_1=&x_5,\quad
z_3=x_4^{-1}y_4,\quad y_3=x_4^2,\\
 z_5=x_6,\quad
&y_5=x_6y_6x_6y_6^{-1}x_6^{-1}y_6^2x_6^2y_6x_6^{-1} >
\end{align*}
Comme nous l'avons dit, $M_1,M_2,M_3,M_4,M_5$ sont des espaces
fibr\'es de Seifert. Dans leur groupes fondamentaux respectifs, la
classe d'une fibre r\'eguli\`ere est respectivement : $x_1$,
$y_2^2$, $y_3$, $x_4^2$ et $y_5$ (ils engendrent le centre).
Remarquons que les isomorphismes $\vf_a$ associ\'es n'envoient pas
la classe d'une fibre r\'eguli\`ere, sur la classe d'une fibre
r\'eguli\`ere, ce qui est en accord avec le fait d\'ej\`a
observ\'e, que la d\'ecomposition $W$ soit minimale.

\section{Extension HNN et amalgame}

Nous donnons dans cette section les th\'eor\`emes de conjugaison
et de commutativit\'e dans un  amalgame ou une extension HNN. Une
preuve dans le cas d'un amalgame peut se trouver dans \cite{mks}.
Nous commen\c cons par rappeler ces résultats dans le cas d'un
amalgame, puis démontrons les résultats similaires dans le cas
d'une extension HNN.

\subsection{Produit amalgam\'e}

Nous consid\'ererons tout au long de ce paragraphe, un groupe
$\G$, produit amalgam\'e des groupes $A$ et $B$ le long du
sous-groupe  $C$.\medskip\\
\textbf{Conjugaison dans un amalgame.} Soit $\w$ un \el\ de $\G$,
donn\'e par une \'ecriture sous forme r\'eduite $\w \equiv \w_1
\w_2\cdots \w_n$. On dira que $\w$ est \textsl{cycliquement
r\'eduit} si  $\w_1$ et $\w_n$ sont dans des facteurs
diff\'erents. Clairement cette d\'efinition ne d\'epend pas du
choix d'une forme r\'eduite de $\w$, ainsi on pourra parler d'\el\
cycliquement r\'eduit.
 De plus, il est clair que tous les conjugu\'es
cycliques $\w_{r}\cdots \w_{n}\w_1 \cdots \w_{r-1}$, d'un \el\
cycliquement r\'eduit, sont aussi cycliquement r\'eduits.

Le r\'esultat suivant constitue le th\'eor\`eme
 4.6 de \cite{mks}.

\begin{thm}
\label{amalgconj1}

Soit $\Gamma = A\, *_{\scriptscriptstyle C} B$. Alors tout
\'el\'ement $\gamma \in \Gamma$ est conjugu\'e \`a un \'el\'ement
cycliquement r\'eduit. De plus si $\gamma$ est un \'el\'ement
cycliquement r\'eduit, alors :
\begin{itemize}
\item[(i)] Si $\gamma$ est conjugu\'e \`a un \'el\'ement $c \in C$,
par $h\in \G$, de forme r\'eduite $h=h_1h_2\ldots h_{n+1}$ avec
$n\geq 0$, alors $\gamma$ est dans un des facteurs, et il existe
une suite finie $(c_0,c_{1},\ldots ,c_{n},\gamma)$, avec $c_0=c$,
et o\`u $\forall\, i=0,\ldots ,n,\ c_{i}\in C$, et $h_{i+1}$
conjugue $c_i$ en $c_{i+1}$ (en posant $c_{n+1}=\g$).
\item[(ii)] Si $\gamma$ est conjugu\'e \`a un \'el\'ement
$\gamma'$, qui est dans un des facteurs, mais n'est pas conjugu\'e
\`a un \'el\'ement de C, alors $\gamma$ et $\gamma'$ sont dans un
m\^eme facteur, et conjugu\'es dans ce facteur.
\item[(iii)] Si $\gamma$ est conjugu\'e \`a un \'el\'ement cycliquement
r\'eduit, d'\'ecriture sous forme r\'eduite $p_{1}p_2\cdots p_{r}$
o\`u $r>1$, alors $\gamma$ peut \^etre obtenu en conjuguant un
conjugu\'e cyclique $p_i\cdots p_rp_1\cdots p_{i-1}$ de
$p_{1}p_2\cdots p_{r}$ par un \'el\'ement de $C$.
\end{itemize}
\end{thm}
On obtient imm\'ediatement  le corollaire suivant :

\begin{cor}[\textbf{Th\'eor\`eme de conjugaison dans un amalgame}]
\label{amalgconj2}  Soient $\g$\linebreak et $\g'$, deux \el\ de
$\G= A\, *_{\scriptscriptstyle C}B$, cycliquement r\'eduits,
conjugu\'es par un \el\ $h$ de $\G$. Alors $\g$ et $\g'$ ont
m\^eme longueur, et de plus :
\begin{itemize}
\item[(i)] Si $\g$ est conjugu\'e \`a un \el\ de C, alors $\g$ et $\g'$ sont de
longueur 1, et si $h$ s'\'ecrit sous forme r\'eduite
$h=h_0h_1\cdots h_{n+1}$, il existe une suite finie
$c_0,c_1,\ldots,c_n$ d'\el s de $C$, telle que, $c_0$ est
conjugu\'e \`a $\g$ par $h_0^{-1}$ dans un facteur, $c_n$ est
conjugu\'e \`a $\g'$ dans un facteur par $h_{n+1}$ , et $\forall\,
i=0,\ldots ,n-1$, $c_i$ et $c_{i+1}$ sont conjugu\'es dans un
facteur par $h_{i+1}$.
\item[(ii)] Si $\g$ est dans un des facteurs et n'est pas conjugu\'e \`a
un \el\ de $C$, alors $\g$ et $\g'$ sont dans le m\^eme facteur,
et conjugu\'es par $h$ dans ce facteur.
\item[(iii)] Si $\g$ n'est pas dans un des facteurs, alors $\g$ s'obtient en conjuguant
un conjugu\'e cyclique de $\g'$ par un \el\ de $C$.
\end{itemize}\smallskip
\end{cor}

\noindent\textbf{Commutativit\'e dans un amalgame.} Le
th\'eor\`eme qui suit, caract\'erise les \'el\'e\-ments qui
commutent dans un produit amalgam\'e. Ses parties (i'), (ii) et
(iii) constituent le th\'eor\`eme 4.5 de \cite{mks}, auquel nous
renvoyons le lecteur pour une preuve.
 Quant \`a la partie (i), elle s'obtient imm\'ediatement en appliquant le
 point (i) du corollaire \ref{amalgconj2}.

\begin{thm}[{\bf Th\'eor\`eme de commutativit\'e dans un amalgame}]
\label{amalgcom}
\hfill Soient $\G=A\, *_{\scriptscriptstyle C} B$, et deux \el s\
$x$ et $y$ de $\G$ qui commutent. Alors :
\begin{itemize}
\item[(i)] Si $x$ est dans $C$, et $y$ s'\'ecrit sous forme r\'eduite,
$y_1 y_2 \cdots y_n$, alors il existe une suite finie
d'\'el\'ements $(c_0,c_1,\cdots ,c_n)$ de $C$, avec $c_0 = c_n =
x$, et
 $\forall\, i=0,\ldots n-1$, $y_{i+1}$ conjugue
$c_i$ en $c_{i+1}$ dans un des facteurs.
\item[(i')] $x$ ou $y$ est dans un conjugu\'e de $C$.
\item[(ii)] Si ni $x$ ni $y$ n'est dans un conjugu\'e de $C$, et
$x$ est dans le conjugu\'e d'un facteur, alors $y$ est dans le
m\^eme conjugu\'e du m\^eme facteur.
\item[(iii)] Si ni $x$ ni $y$ n'est dans le conjugu\'e d'un facteur,
alors $x=g h g^{-1} .W^j$, et $y= g h' g^{-1} .W^k$, o\`u $g,W \in
\G$, et $h,h' \in C$, et $ghg^{-1} ,gh'g^{-1}$, et $W$, commutent
deux \`a deux.
\end{itemize}
\end{thm}

Quant au centre d'un amalgame il est caractérisé par le résultat
suivant (cor.4.5, \cite{mks}) :

\begin{thm}[{\bf Centre d'un amalgame}]
\label{amalgcenter}
\hfill Le centre d'un amalgame $\G=A\, *_{\scriptscriptstyle C} B$
non trivial (\textsl{i.e.} $A\not=C$, $B\not=C$)  est
$Z(G)=Z(A)\cap Z(B)$.
\end{thm}

\subsection{Conjugaison dans une extension HNN}
Soit  le groupe $\G$, extension HNN de $A$ le long de
l'isomorphisme $\f :C_{-1}\longrightarrow C_{+1}$. Soit $\w$ un
\el\ de $\G$, qui s'\'ecrit sous forme r\'eduite
 $\w\equiv \w_1t^{\e_1}\cdots t^{\e_n}\w_{n+1}$.
On dit que $\w$ est \textsl{cycliquement r\'eduit}
 si
soit $|\w |=1$,  soit $\w_{n+1}=1$, et dans ce dernier cas, soit
$|\w |=2$, (\emph{i.e.} $w=w_1t^{\e_1}$), soit $|\w |>2$, et
$t^{\e_{n}}\w_{1}t^{\e_{1}}$ n'est pas un pinch. Il est clair que
cette condition ne d\'epend pas du choix d'une forme normale
r\'eduite de $\w$. Aussi on pourra parler d'\el\ cycliquement
r\'eduit de $\G$.

Le th\'eor\`eme de conjugaison que nous donnons n'est qu'une forme
plus explicite du lemme de Collins (th\'eor\`eme 2.5, chap. IV
\cite{Lyndon}) ; aussi nous n'en fournirons une preuve qu'en
appendice (\S 7).

\begin{thm}
\label{hnnconj1}
Soit $\Gamma=A\ast_\f$ avec $\f:C_{-1}\longrightarrow C_{+1}$ ;
alors tout \'el\'ement $\g\in \G$ est conjugu\'e \`a un
\'el\'ement cycliquement r\'eduit. De plus si $\g$ est
cycliquement r\'eduit, alors :

\begin{itemize}

\item[(i)] Si $\g$ est conjugu\'e \`a un \'el\'ement
$c\in C_{+1}\cup C_{-1}$ par l'\el\ $h$ de forme r\'eduite
$h=h_1t^{\e_1}h_2\cdots t^{\e_p}h_{p+1}$,
 alors $\g$ est dans $A$,
et il existe une suite finie $(c_{0},c_{1},\cdots,c_{2p})$
d'\'el\'ements de $C_{+1}\cup C_{-1}$,
 telle que $c_{2p}=c$,
$\g=h_1.c_0.h_1^{-1}$, et pour $i=0,\ldots ,p-1$,
\begin{gather*}
c_{2i+1}=h_{i+2}. c_{2i+2}. h_{i+2}^{-1} \\
c_{2i}= t^{\e_{(i+1)}}. c_{2i+1}. t^{-\e_{(i+1)}}
\end{gather*}

\item[(ii)] Si $\g$ est conjugu\'e \`a un \'el\'ement $\g'\in A$, mais
n'est pas conjugu\'e \`a un \'el\'ement de $C_{+1}\cup C_{-1}$,
alors  $\g$ est dans $A$, et $\g$ et $\g'$ sont conjugu\'es dans
$A$.

\item[(iii)] Si $\g$ est conjugu\'e \`a un \'el\'ement $\g'$,
cycliquement r\'eduit, d'\'ecriture sous forme r\'eduite,
$u_{1}t^{\mu_{1}} \cdots u_{m}t^{\mu_{m}}$ avec $m\geq 1$, alors
$\g$ peut-\^etre obtenu en conjuguant un conjugu\'e cyclique de
$\g'$, $u_{r+1}t^{\mu_{r+1}}\cdots u_mt^{\mu_m}u_1t^{\mu_1}\cdots
u_{r} t^{\mu_{r}}$ par un \'el\'ement de $C_{\mu_{r}}$.

\end{itemize}
\end{thm}
Et l'on obtient alors imm\'ediatement le corollaire suivant :

\begin{cor}[{\bf Th\'eor\`eme de conjugaison dans une extension HNN}]
\label{hnnconj2}
Soit $\Gamma=A\ast_\f$ avec $\f:C_{-1}\longrightarrow C_{+1}$ ; si
$\g$ et $\g'$ sont deux \el s cycliquement r\'eduits de $\G$
 conjugu\'es par $h\in \G$, alors $\g$ et $\g'$
ont m\^eme longueur, et de plus :
\begin{itemize}
\item[(i)] Si $\g$ est conjugu\'e \`a un \el\ $c$ de $C_{+1} \cup C_{-1}$,
alors $\g$ et $\g'$ sont dans $A$, et si $h$ s'\'ecrit sous la
forme r\'eduite  $h=h_1t^{\e_1}\cdots t^{\e_p}h_{p+1}$, il existe
une suite finie $(c_0,c_1,\ldots,c_{2p-1})$ d'\el s de $C_{+1}
\cup C_{-1}$, telle que $\g=h_1.c_0.h_1^{-1}$, et pour $i=0,\ldots
,p-1$, en posant $c_{2p}=\g'$,
\begin{gather*}
c_{2i+1}=h_{i+2}. c_{2i+2}. h_{i+2}^{-1} \\
c_{2i}= t^{\e_{(i+1)}}. c_{2i+1}. t^{-\e_{(i+1)}}
\end{gather*}

\item[(ii)] Si $\g$ est dans $A$, mais n'est pas conjugu\'e \`a un
\el\ de $C_{+1} \cup C_{-1}$, alors $\g'$ est dans $A$ et $\g$ et
$\g'$ sont conjugu\'es par $h$ dans $A$.
\item[(iii)] Si $\g$ n'est pas dans $A$,
$\g = u_1 t^{\mu_1}\cdots u_m t^{\mu_m}$, avec $m\geq 1$ alors
$\g$ s'obtient en conjuguant un conjugu\'e cyclique
$v_{r}t^{\mu_{r}}\cdots v_m t^m v_1 t^1 \cdots v_{r-1}
t^{\mu_{r-1}}$ de $\g'$ par un \el\ $\a$ de $C_{\mu_m}$. \medskip
\end{itemize}
\end{cor}

\subsection{Commutativit\'e dans une extension HNN}
Le théorème suivant caractérise deux éléments qui commutent dans
une extension HNN.

\begin{thm}[{\bf Th\'eor\`eme de commutativit\'e dans une extension HNN}]
\label{hnncom}  Soient $\G = A\ast_\f$ avec $\f : C_{-1}
\longrightarrow C_{+1}$, et $x$ et $y$, deux \el s\ de $G$ qui
commutent. Alors :
\begin{itemize}
\item[(i)] Si $x\in C_{+1} \cup C_{-1}$, et
$y$ s'\'ecrit sous forme r\'eduite, $y=y_1 t^{\e_1} y_2 \cdots y_n
t^{\e_n} y_{n+1}$, alors il existe une suite d'\el s\  $(c_0 ,c_1,
\cdots , c_{2n+1})$ de $C_{+1} \cup C_{-1}$ o\`u $c_0 = c_{2n+1} =
x$, et $c_i$ et $c_{i+1}$, sont conjugu\'es, par
$y_{\frac{i}{2}+1}$ si $i$ est pair, et par $t^{\e_j}$ avec $j=
{(\frac{i-1}{2}+1)}$, si $i$ est impair.
\item[(i')] x ou y est conjugu\'e \`a un \'el\'ement de $C_{+1}\cup C_{-1}$.
\item[(ii)] Si $x$ est dans un conjugu\'e de $A$ et n'est pas conjugu\'e \`a
un \el\ de $C_{+1} \cup C_{-1}$, alors $y$ est dans le m\^eme
conjugu\'e de $A$.
\item[(iii)] Si ni $x$ ni $y$ n'est dans un conjugu\'e de $A$,
alors $x=gcg^{-1} W^j$ et $y=gc'g^{-1}W^k$, o\`u $g,W\in \G$,
$c,c'\in C_{+1}$ (respectivement $c,c' \in C_{-1}$), et
$gcg^{-1}$, $gc'g^{-1}$, et $W$ commutent deux \`a deux.
\end{itemize}
\end{thm}

\noindent \textbf{Démonstration.} Nous traitons sépar\'ement
chacun des
cas.\\
\textsl{Cas} (i). Il suffit d'appliquer le corollaire
\ref{hnnconj2}, apr\`es avoir remarqu\'e que $y$ conjugue $x$ en
$x$.
\hfill $\square$\smallskip\\
\noindent \textsl{Cas} (i'). Il n'y a rien \`a prouver dans ce
cas.
\hfill $\square$\smallskip\\
\noindent \textsl{Cas} (ii). On suppose que $x$ est dans $A$. On a
encore
\begin{align*}
x &= y \; x \; y^{-1}  \\
x &= y_1 t^{\e_1} y_2 \cdots y_n t^{\e_n} y_{n+1} \; x \;
y_{n+1}^{-1} t^{-\e_n}
 \cdots y_2 ^{-1} t^{-\e_1} y_1^{-1}
\end{align*}
qui est dans $A$, et donc le membre de droite n'est pas r\'eduit.
Or puisque $x$ n'est pas conjugu\'e \`a un \el\ de $C_{+1} \cup
C_{-1}$, $t^{\e_n} y_{n+1} \; x \; y_{n+1}^{-1} t^{-\e_n}$ ne peut
pas \^etre un pinch, et donc, puisque $y$ est r\'eduit, la seule
possibilit\'e est $n=0$, c'est \`a dire
$y \in A$.\\
Si $x$ est conjugu\'e \`a un \el\ $x_0 \in A$, $x= h \; x_0 \;
h^{-1}$, alors
$h^{-1} \; y \; h$ commute avec $x_0$, et est donc dans $A$.\hfill $\square$\smallskip\\
\indent Pour achever la d\'emonstration, nous proc\'edons par
contradiction. Nous supposons qu'il existe $x,y\in \G$, qui
commutent, et qui ne v\'erifient pas les conclusions de la
proposition. Alors, avec tout ce qui pr\'ec\`ede, ni $x$ ni $y$
n'est dans un conjugu\'e de $A$. En particulier, $x$ et $y$ sont
de longueur sup\'erieure \`a 1.

On consid\`ere un \'el\'ement $x\in \G$, de longueur minimale,
pour lequel il existe $y\in \G$, tel que $x$ et $y$ commutent et
ne v\'erifient pas les conclusions de la proposition. On prend $y$
de longueur minimale dans l'ensemble des \'el\'ements v\'erifiant
cette condition. On pose $r+1=|x|,s+1=|y|$ ; $x$ et $y$
s'\'ecrivent sous forme r\'eduite :
$$x=x_1t^{\e_1}x_2\cdots x_rt^{\e_r}x_{r+1}$$
$$y=y_1t^{\mu_1}y_2\cdots y_st^{\mu_s}y_{s+1}$$
Bien s\^ ur, $0<r\leq s$.

Remarquons tout d'abord, que $t^{\e_r}x_{r+1}x_1t^{\e_1}$ n'est
pas un pinch. En effet, dans le cas contraire, $t^{\e_r}x_{r+1}\;
x\; x_{r+1}^{-1}t^{-\e_r}$ est de longueur strictement
inf\'erieure \`a $|x|$, et $t^{\e_r}x_{r+1}\; x\;
x_{r+1}^{-1}t^{-\e_r}$ et $t^{\e_r}x_{r+1}\; y\;
x_{r+1}^{-1}t^{-\e_r}$ commutent, et ne v\'erifient pas les
conclusions de la proposition, ce qui est contradictoire, puisque
l'on a suppos\'e, que $x$ \'etait de longueur minimale.

Maintenant, soit $xy$, soit $xy^{-1}$ est sous forme r\'eduite de
longueur $r+s+1$. Dans le cas contraire,
$t^{\e_r}x_{r+1}y_1t^{\mu_1}$ et $t^{-\e_1}x_1^{-1}y_1t^{\mu_1}$
sont des pinchs. Alors, $\e_r=-\e_1$, et $x_{r+1}y_1,
x_1^{-1}y_1\in C_{\e_r}$, et donc :
$$x_{r+1}y_1=c \; x_1^{-1}y_1$$
o\`u $c\in C_{\e_r}$, et alors
$$x_{r+1}x_1=c$$
et $t^{\e_r}x_{r+1}x_1t^{\e_1}$ est un pinch, ce qui est
contradictoire.

Sans perte de g\'en\'eralit\'e, on peut supposer que $xy$ est sous
forme r\'eduite de longueur $r+s+1$. Ainsi,
\begin{align*}
xy&= x_1t^{\e_1}\cdots t^{\e_r}x_{r+1}\, y_1t^{\mu_1}
\cdots t^{\mu_s}y_{s+1}\qquad\qquad (*)\\
&= yx\\
&=y_1t^{\mu_1}\cdots t^{\mu_s}y_{s+1}\, x_1t^{\e_1} \cdots
t^{\e_r}x_{r+1} \qquad\qquad (**)
\end{align*}
et $(*)$ et $(**)$ sont deux \'ecritures r\'eduites de $xy=yx$.
Avec le lemme de Britton, en supposant $s\geq r$,
$$t^{\mu_{q+1}}y_{q+2}\cdots t^{\mu_s}y_{s+1}=c\;t^{\e_1}x_2\cdots t^{\e_r}x_{r+1}$$
o\`u $q=s-r$, et $c\in C_{-\e_1}$. Alors,
$$yx^{-1}=y_1t^{\mu_1}\cdots t^{\mu_q}y_{q+1}\, c \, x_1^{-1}$$
est de longueur $q+1=s-r+1$. Or, $x$ et $yx^{-1}$ commutent, et
donc puisque $|yx^{-1}|<|y|$, avec le choix que nous avons fait de
$x$ et $y$, n\'ecessairement $x$ et $yx^{-1}$ v\'erifient les
conclusions de la proposition.

Si $yx^{-1}$ est conjugu\'e \`a un \'el\'ement $c\in C_{+1}\cup
C_{-1}$,
$$yx^{-1}=gcg^{-1}$$
et donc,
$$y=gcg^{-1}\, x$$
et puisque $x=g1g^{-1}\, x$, $x$ et $y$ v\'erifient (iii), ce qui
est contradictoire.

Si $yx^{-1}$ est dans un conjugu\'e de $A$, et n'est pas
conjugu\'e \`a un \el\ de $C_{+1}\cup C_{-1}$, avec (ii), $x$ et
$y$ sont dans le m\^eme conjugu\'e de $A$, et v\'erifient donc
(ii), ce qui est contradictoire.

Si $yx^{-1}=g c g^{-1}\, W^j$, et $x= g c' g^{-1}\, W^k$,
v\'erifient la conclusion (iii), alors
\begin{align*}
y=yx^{-1}\, x &= g c g^{-1} \, W^j \; g c' g^{-1} \, W^k\\
    &= g c g^{-1} \, g c' g^{-1} \, W^j\, W^k\\
    &= g cc' g^{-1} \, W^{j+k}
\end{align*}
et puisque $g c g^{-1}, g c' g^{-1}, W$ commutent deux \`a deux,
$g c' g^{-1}, g cc' g^{-1}$ et $W$ commutent deux \`a deux. De
plus $c,c'$ sont dans $C_{+1}$ (respectivement $C_{-1}$), et donc
$cc'$ est dans $C_{+1}$ (respectivement $C_{-1}$). Ainsi, $x$ et
$y$ v\'erifient la conclusion (iii), ce qui est contradictoire. Il
n'existe donc pas d'\el s $x,y\in \G$, qui commutent et ne
v\'erifient pas les conclusions de la proposition. Ceci termine la
d\'emonstration.
 \hfill $\blacksquare$
\vskip 0.4cm

\subsection{Centre d'une extension HNN}
Le résultat suivant caractérise le centre d'une extension HNN :

\begin{thm}[{\bf Centre d'une extension HNN}]
\label{hnncenter}
 Soit $\G = A\ast_\f$, avec $\f : C_{-1}
\longrightarrow C_{+1}$ ; notons $\mathrm{Fix}\ \f$ le sous-groupe
de $C_{-1}\cap C_{+1}$ constitué des point fixes de $\f$, et
$Z(\G)$ le centre de $\G$.
\begin{itemize}
\item[(i)] Si $C_{-1}$ ou $C_{+1}$ est un sous-groupe propre de
$A$, alors $Z(\G)=Z(A)\cap \mathrm{Fix}\ \f$.
\item[(ii)] Si $C_{-1}=C_{+1}=A$ et si $\f$ est d'ordre infini
dans $Out(A)$, alors $Z(\G)=Z(A)\cap \mathrm{Fix}\ \f$.
\item[(iii)] Si $C_{-1}=C_{+1}=A$ et si $\f$ est d'ordre $n$ dans
$Out(A)$ alors soit $a_0\in \G$ tel que pour tout $a_0\in A$,
$\f^n(a_0)=a_0aa_0^{-1}$. Alors $Z(\G)=\bigcup_{p\in \Z}
t^{pn}.(\mathrm{Fix}\ \f\cap a_0^pZ(A))$.
\end{itemize}
\end{thm}

\noindent{\bf Démonstration.} Remarquons au préalable que puisque
$Z(A)\cap \mathrm{Fix}\ \f$ est dans le centralisateur de la
famille g\'en\'eratrice $A\cup\{ t\}$ de $\G$ on a dans tous les
cas $Z(A)\cap\mathrm{Fix}\ \f\subset Z(\G)$. Soit $x\in Z(\G)\cap
A$ ; nécessairement, d'une part $x\in Z(A)$, et d'autre part
$t^{-1}xt=x$, et donc $x=\f(x)$. Ainsi dans tous les cas
$Z(\G)\cap
A = Z(A)\cap\mathrm{Fix}\ \f$.\smallskip\\
\indent Nous commençons par traiter le cas (i) : $C_{-1}$ ou
$C_{+1}$ est propre dans $A$. Puisque l'isomorphisme naturel entre
$A*_\f$ et $A*_{\f^{-1}}$ pr\'eserve $A$ et que $\mathrm{Fix}\
\f=\mathrm{Fix}\ \f^{-1}$ on se restreindra sans perte de
g\'en\'eralit\'e au cas où $C_{-1}$ est propre dans $A$. Avec ce
qui précède il nous suffit de montrer que $Z(\G)\subset A$.

 Soit $x\in Z(\G)$ ; ab
absurdo supposons que $|x|>1$. Remarquons tout d'abord que puisque
$x$ est dans le centre de $\G$,  $x$ \'egale tous ses conjugu\'es
cycliques. Ainsi d'une part $x$ admet une écriture cycliquement
réduite :
$$x=x_1t^{\e_1}\ldots x_nt^{\e_n}$$
avec $n\geq 1$, et d'autre part tous les $\e_i$ sont \'egaux. Considérons $y\in A\setminus C_{-1}$.\\
Si $\e_n=-1$, alors $xyx^{-1}$ a une écriture réduite de longueur
$2n+1$, et est égal à $y$ qui est de longueur 1 ; ceci est
contraditoire.\\
Si $\e_n=\e_1=1$,  on applique le m\^eme raisonnement que
ci-dessus en consid\'erant au lieu de $x$,
$x^{-1}=x_1^{-1}t^{-\e_n}x_n^{-1}\ldots x_2^{-1}t^{-\e_1}$, pour
arriver encore à une contradiction ; ceci montre que
n\'ecessairement $x\in A$, ce qui achève la preuve du cas
(i).\hfill $\square$

Nous traitons maintenant les deux derniers cas ; remarquons que
lorsque $C_{-1}=C_{+1}=A$, $\G$ est le produit semi-direct
$\G=A\rtimes_\f \Z$, et que tout \el\ de $\G$ s'\'ecrit de fa\c
con unique $t^na$ avec $a\in A$ et $n\in\Z$.

Si $x=t^na_1\in Z(\G)$, alors  pour tout $a\in A$, $a
t^na_1a^{-1}=t^n\f^n(a)a_1a^{-1}=t^na_1$, et donc $\f^n(a)=a_1a
a_1^{-1}$. Ainsi soit $Z(\G)\subset A$, soit $\f$ est d'ordre fini
dans $Out(A)$. En particulier dans le cas (ii) on a la conclusion
souhaitée. Il ne nous reste donc plus qu'à montrer le cas (iii).

Supposons que $\f$ est d'ordre $n$ dans $Out(A)$ et considérons
$a_0\in A$ tel que pour tout $a\in A$, $\f^n(a)=a_0a a_0^{-1}$ ;

 Soient $p\in
\Z$, $\a\in Z(A)$ tel que $a_0^p\a \in \mathrm{Fix}\ \f$, et
$x=t^{pn}a_0^p\a$ ; soit $y$ un \el\ quelconque de $\G$, $y=t^ra$
avec $y\in \Z$ et $a\in A$.
\begin{align*}
yxy^{-1} &= t^ra.t^{pn}a_0^p\a.a^{-1}t^{-r}\\
&= t^{pn+r}a_0^paa_0^{-p}.a_0^{p}\a a^{-1}t^{-r}\\
&=t^{pn+r}a_0^p\a t^{-r}\\
&=t^{pn}a_0^p\a=x\\
\end{align*}
 Ceci montre que $\bigcup_{p\in\Z}
t^{pn}.(\mathrm{Fix}\ \f\cap a_0^pZ(A))\subset Z(\G)$.

R\'eciproquement soient $m\in \Z$, $a_1\in A$ et $x=t^ma_1\in
Z(G)$. Alors d'une part $\forall a\in A$, $\f^m(a)=a_1aa_1^{-1}$
ce qui implique qu'il existe $p\in\Z$ tel que $m=np$ et $a_1\in
a_0^pZ(A)$. D'autre part $t^{-1}xt=x$ montre que $t^m\f
(a_1)=t^ma_1$ et donc que $a_1\in \mathrm{Fix}\ \f$. Ainsi
 $Z(\G)\subset \bigcup_{p\in\Z}t^{pn}.(\mathrm{Fix}\ \f\cap a_0^pZ(A))$ ce
qui achève la preuve de (iii).\hfill $\blacksquare$


\section{Graphes de groupe}

Nous commençons par \'etablir dans cette section  les
th\'eor\`emes principaux caractérisant le centre, les éléments qui
commutent et les éléments conjugués
 (théorèmes \ref{graph center}, \ref{graph_com} et \ref{graph_conj}) dans le groupe fondamental d'un graphe de groupe.

L'approche employée est inductive. Si $(\cal{G},X)$ est un graphe
 de groupe, d\'ecomposer $X$ le long d'une de ses ar\^etes,
d\'ecompose $\P(\cal{G},X)$ en extension HNN ou en amalgame de
groupes fondamentaux de sous-graphes(s) de groupe. Les
th\'eor\`emes de la section pr\'ec\'edente servent alors \`a la
fois de conditions initiales, et de conditions de r\'ecurrence.

 Nous  établirons ensuite divers résultats. D'abord nous montrons  que
pour un groupe quelconque $G$  et une  famille finie quelconque de
ses sous-groupes, le problème de conjugaison dans $G$ se réduit au
problème de conjugaison dans le double $2G$ de $G$ le long de ces
sous-groupes (théorème \ref{double}). Ensuite nous donnons une
condition combinatoire suffisante  pour que dans un groupe
fondamental de graphe de groupe, les centralisateurs, le centre et
la structure des racines soient dans un sens triviaux (théorèmes
\ref{centralizer}, \ref{triv center} et \ref{srt}).

\subsection{Centre dans un graphe de groupe}

Nous donnons d'abord quelques définitions nécessaires à l'énoncé
du théorème.

\begin{defn}
Un graphe de groupe décomposé $(\cal G, X, T)$ est dit {\sl
minimal} si pour toute arête $T$-séparante $a$, le sous-groupe
d'arête $G_a$ est un sous-groupe propre des sous-groupes de
sommets $G_{\mathrm o(a)}$ et $G_{\mathrm e(a)}$ de $\pi_1(\cal
G,X)$.
\end{defn}

\noindent{\bf Remarque :} Donné un graphe de groupe décomposé qui
n'est pas minimal on peut clairement --en 'écrasant' en des
sommets les arêtes qui posent problème-- construire un graphe de
groupe décomposé minimal, sans en changer le groupe fondamental.

\begin{defn}
Soit un graphe de groupe décomposé $(\cal G, X, T)$ ; donnée une
arête $a$ non $T$-séparante, on note $\mathrm{Fix}\ a$ le
sous-groupe de $G_a^-\cap G_a^+\subset \pi_1(\cal G,X)$ constitué
des éléments fixés par l'isomorphisme $\vf_a :
G_a^-\longrightarrow G_a^+$. (Remarquons que clairement
$\mathrm{Fix}\ a=\mathrm{Fix}\ {-a}$).
\end{defn}

Nous pouvons dès-lors donner le théorème caractérisant le centre
du groupe fondamental d'un graphe de groupe.

\begin{thm}[{\bf Centre dans un graphe de groupe}]
\label{graph center}
 Soient $(\cal G,X,T)$ un graphe de groupe
décomposé minimal et $\G$ son groupe fondamental. Soit $X$ a un
unique sommet et une unique arête non-orientée (dans ce cas la
conclusion est donnée par le théorème \ref{hnncenter}), soit :
$$Z(\G)=(\bigcap_{s\in\cal S_X} Z(G_s))\cap (\bigcap _{a\in\cal
A_{X\setminus T}} \mathrm{Fix}\ a)$$
\end{thm}

\noindent{\bf Démonstration.} Nous procédons par induction sur le
nombre d'arêtes non $T$-séparantes. Si $X$ n'a pas d'arête non
$T$-séparante (\emph{i.e.} $X=T$), alors avec le théorème
\ref{amalgcenter}, $Z(G)=\bigcap_{s\in\cal S_X} Z(G_s)$. Supposons
maintenant que $X$ contienne au moins une arête non $T$-séparante
$a$, et considérons le sous-graphe de groupe décomposé $(\cal
G,X',T)$ où $X'$ est obtenu à partir de $X$ en supprimant l'arête
$a$. Alors $\G$ est l'extension HNN de base $\G_0=\pi_1(\cal
G,X')$ le long de l'isomorphisme $\vf_a : G_a^-\longrightarrow
G_a^+$. Si $X$ n'a pas un unique sommet et une unique arête
non-orientée, alors avec l'hypothèse de minimalité $G_a^-$ est un
sous-groupe propre de $\G_0$, ainsi le théorème \ref{hnncenter}
(i) s'applique pour montrer que $Z(\G)=Z(\G_0) \cap \mathrm{Fix}\
a$, ce qui achève l'induction.\hfill $\blacksquare$

\subsection{Trajets et circuits} Nous introduisons dans cette
section le matériel combinatoire qui nous sera utile dans la
suite.\\

\begin{defn}
 Soit $(\cal{G},X,T)$ un
graphe de groupe d\'ecompos\'e. Soient $s_o$ et $s_e$ des sommets
de $X$,
 et les
\el s  $u\in G_{s_o}$, $v\in G_{s_e}$ de $\P(\cal{G},X)$. Un
\textsl{trajet} de $u$ \`a $v$, d'origine $s_o$
et d'extr\'emit\'e $s_e$, est la donn\'ee de :\medskip\\
\indent a) Un chemin de $X$, d'origine $s_o$ et d'extr\'emit\'e
$s_e$,
$$(a_1,a_2,\ldots ,a_n)$$
avec $n\geq 0$ (si $n=0$ le chemin est réduit à $s_0$). On note
$s_0=s_o, s_{n}=s_e$, et si $n\geq 2$, pour tout $i=1,\ldots
,n-1$,
$s_i=\mathrm{e}(a_i)$.\medskip\\
\indent b) Une suite
$$(c_1^-,c_1^+,c_2^-,c_2^+,\ldots ,c_n^-,c_n^+)$$
avec $\forall\, i=1,\ldots ,n$, $c_i^-\in G_{a_i}^-$,
$c_i^+\in G_{a_i}^+$, et $c_i^+=\varphi_{a_i}(c_i^-)$.\medskip\\
\indent c) Une suite
$$(h_0,h_1,\ldots h_{n})$$
avec soit $n=0$ et $u=h_0vh_0^{-1}$ dans $G_{s_0}$ ; soit
 $\forall\, i= 0,\ldots ,n$, $h_i\in G_{s_i}$, et
$$u=h_0c_1^-h_0^{-1}\quad \mathrm{dans\;} G_{s_0}$$
$$c_n^+= h_nvh_n^{-1}\quad \mathrm{dans\;} G_{s_n}$$
et si $n\geq 2$, $\forall\, i=1,\ldots , n-1$,
$$c_i^+=h_ic_{i+1}^-h_i^{-1}\quad \mathrm{dans\;} G_{s_i}$$
On symbolisera la donn\'ee d'un tel trajet, par la notation
suivante :
$$u\underset{h_0}{\circlearrowleft} c_1^-
\overset{a_1}{\longrightarrow} c_1^+
\underset{h_1}{\circlearrowleft} c_2^-
\overset{a_2}{\longrightarrow} \cdots
\overset{a_i}{\longrightarrow} c_{i}^+
\underset{h_i}{\circlearrowleft} c_{i+1}^-
\overset{a_{i+1}}{\longrightarrow} \cdots
\overset{a_{n-1}}{\longrightarrow}
c_{n-1}^+\underset{h_{n-1}}{\circlearrowleft}c_n^-
\overset{a_n}{\longrightarrow} c_n^+
\underset{h_n}{\circlearrowleft} v$$
Si $n=0$, on parlera du trajet  \textsl{trivial},
 que l'on notera :
$$u\underset{h_0}{\circlearrowleft} v$$
Si $h_i\in G_{s_i}$ est l'\el\ neutre, on notera
$$\cdots c_i^+= c_{i+1}^- \cdots$$
au lieu de
$$\cdots c_i^+ \underset{1}{\circlearrowleft} c_{i+1}^- \cdots$$
\end{defn}

\begin{defn}\label{traj} Si le trajet $\cal{C}$ est donn\'e par :
$$u\underset{h_0}{\circlearrowleft} c_1^-
\overset{a_1}{\longrightarrow} c_1^+
\underset{h_1}{\circlearrowleft} c_2^-
\overset{a_2}{\longrightarrow} \cdots
\overset{a_i}{\longrightarrow} c_{i}^+
\underset{h_i}{\circlearrowleft} c_{i+1}^-
\overset{a_{i+1}}{\longrightarrow} \cdots
\overset{a_{n-1}}{\longrightarrow}
c_{n-1}^+\underset{h_{n-1}}{\circlearrowleft}c_n^-
\overset{a_n}{\longrightarrow} c_n^+
\underset{h_n}{\circlearrowleft} v$$ on appelle \textsl{label} du
trajet $\cal{C}$,  l'\el\ $h\in\P(\cal{G},X)$, que l'on notera
$\mathrm{label}(\cal{C})$, d\'efini par
$$\mathrm{label}(\cal{C})=h_0t_{a_1}h_1t_{a_2}\cdots t_{a_n}h_n=h$$\label{label}
\noindent pour la pr\'esentation donn\'ee par le th\'eor\`eme
\ref{graphpres} (rappelons que si $a$ est une ar\^ete
$T$-s\'eparante, $t_a=t_{-a}=1$). Il est clair, qu'alors
$$u=hvh^{-1}\quad \mathrm{dans\; }\P(\cal{G},X)$$
\end{defn}

\begin{defn}
Un \textsl{circuit} en $u$
d'origine $s_o$,  est un trajet pour lequel $s_o=s_e$, de $u$ \`a
$u$. Si $h$ est l'\el\ associ\'e \`a un circuit en $u$, alors
$$\lbrack u,h\rbrack =1\quad \mathrm{dans\; } \P(\cal{G},X)$$
Un \textsl{circuit trivial} en $u$ est un circuit de la forme :
$$u\underset{h_0}{\circlearrowleft}u$$
\end{defn}
\begin{defn}
Consid\'erons un trajet
$\cal{T}_1$ de $u$ \`a $v$, d'origine $s_1$ et d'extr\'emit\'e
$s_2$.
$$u\underset{h_0}{\circlearrowleft} c_1^-
\overset{a_1}{\longrightarrow} \cdots
\overset{a_i}{\longrightarrow} c_{i}^+
\underset{h_i}{\circlearrowleft} c_{i+1}^-
\overset{a_{i+1}}{\longrightarrow} \cdots
\overset{a_n}{\longrightarrow} c_n^+
\underset{h_n}{\circlearrowleft} v$$ Notons $\cal{T}_1^{-1}$ le
trajet ainsi d\'efini :
$$v\underset{h_n^{-1}}{\circlearrowleft} c_n^+
\overset{-a_n}{\longrightarrow} \cdots
\overset{-a_{i+1}}{\longrightarrow} c_{i+1}^-
\underset{h_i^{-1}}{\circlearrowleft} c_{i}^+
\overset{-a_i}{\longrightarrow} \cdots
\overset{-a_1}{\longrightarrow} c_1^-
\underset{h_0^{-1}}{\circlearrowleft} u$$ C'est un trajet de $v$
\`a $u$, d'origine $s_2$ et d'extr\'emit\'e $s_1$. Si $h_1$ est le
label de $\cal{T}_1$, alors $\cal{T}_1^{-1}$ a pour label
$h_1^{-1}$. On dira que le trajet $\cal{T}_1^{-1}$ est
\textsl{l'inverse} de $\cal{T}_1$.
 Consid\'erons un trajet $\cal{T}_2$ d'origine $s_2$ et
d'extr\'emit\'e $s_3$, de $v$ \`a $w$ :
$$v\underset{k_0}{\circlearrowleft} d_1^-
\overset{b_1}{\longrightarrow} \cdots
\overset{b_j}{\longrightarrow} d_{j}^+
\underset{k_j}{\circlearrowleft} d_{j+1}^-
\overset{b_{j+1}}{\longrightarrow} \cdots
\overset{b_m}{\longrightarrow} d_m^+
\underset{k_m}{\circlearrowleft} w$$ Sous ces hypoth\`eses, on
d\'efinit le \textsl{produit} de $\cal{T}_1$ et $\cal{T}_2$,
 not\'e $\cal{T}_1\cal{T}_2$, par
$$u\underset{h_0}{\circlearrowleft} c_1^-
\overset{a_1}{\longrightarrow} \cdots
\overset{a_i}{\longrightarrow} \cdots
\overset{a_n}{\longrightarrow} c_n^+
\underset{h_nk_0}{\circlearrowleft} d_1^-
\overset{b_1}{\longrightarrow} \cdots
\overset{b_j}{\longrightarrow} \cdots
\overset{b_m}{\longrightarrow} d_m^+
\underset{k_m}{\circlearrowleft} w$$ C'est un trajet de $u$ \`a
$w$, d'origine $s_1$ et d'extr\'emit\'e $s_3$. Si $h_2$ est le
label de $\cal{T}_2$, alors le label de $\cal{T}_1\cal{T}_2$ est
le produit $h_1h_2$. Il est facile  de v\'erifier que lorsqu'il
est d\'efini, ce produit est associatif.
\end{defn}

\begin{defn}
\label{C(u,s)} 
Consid\'erons un \el\ non trivial $u$ dans un sous-groupe de
sommet $G_s$, et notons $\cal{C}(u,s)$ l'ensemble des labels des
circuits en $u$ d'origine $s$. Remarquons que le produit de deux
circuits en $u$ d'origine $s$ est toujours d\'efini, et que le
produit et l'inverse de deux circuits en $u$ d'origine $s$, est
encore un circuit en $u$ d'origine $s$. Ainsi $\cal{C}(u,s)$ muni
de l'op\'eration de $\P(\cal{G},X)$, forme un sous-groupe de
$\P(\cal{G},X)$.
\end{defn}

Les premi\`eres propositions que nous d\'emontrons, ont pour but
de ne plus avoir \`a parler de l'origine et de l'extr\'emit\'e
d'un trajet.

\begin{prop}
\label{III1} Soit $(\cal{G},X,T)$ un graphe de groupe
d\'ecompos\'e.
 Soient $s_1,s_2$ deux sommets distincts
de $X$, et $u$ un \el\  de $\P(\cal{G},X)$.

Alors $u\in G_{s_1}\cap G_{s_2}$, si et seulement si il existe un
trajet dans $\P(\cal{G},X)$, de $u$ \`a $u$, d'origine $s_1$, et
d'extr\'emit\'e $s_2$, de la forme
$$u = c_1^-
\overset{a_1}{\longrightarrow} \cdots
\overset{a_i}{\longrightarrow} c_{i}^+ = c_{i+1}^-
\overset{a_{i+1}}{\longrightarrow} \cdots
\overset{a_n}{\longrightarrow} c_n^+ = u$$ avec pour tout
$i=1,\ldots ,n$, $a_i\in \cal{A}(T)$ (et donc le label
 de ce trajet est l'\el\ neutre).
\end{prop}

\noindent\textbf{D\'emonstration.} La r\'eciproque \'etant
triviale, nous d\'emontrons l'implication directe.

Remarquons qu'un trajet dans un sous-graphe de groupe de
$(\cal{G},X)$ est aussi naturellement un trajet de $(\cal{G},X)$,
ayant m\^eme label. Ainsi, puisque $G_{s_1},G_{s_2}\subset
\P(\cal{G},T)\subset\P(\cal{G},X)$, il est suffisant de montrer la
propri\'et\'e dans le sous-graphe $(\cal{G},T)$ de $(\cal{G},X)$,
muni de la d\'ecomposition induite.  Consid\'erons dans $T$
l'unique chemin r\'eduit $(a_1,\ldots ,a_n)$ d'origine $s_1$ et
d'extr\'emit\'e $s_2$, et notons $C$ le sous-graphe orient\'e de
$T$ d'ar\^etes $a_1,\ldots, a_n$. Comme pr\'ec\'edemment,
$\P(\cal{G},C)$ contient $G_{s_1}$ et $G_{s_2}$, et se plonge
naturellement dans $\P(\cal{G},T)$, aussi on travaille dans
$(\cal{G},C)$.

Soit $u\in G_{s_1}\cap G_{s_2}$. D\'ecomposons $C$ le long de
$a_1$. Notons $C_1$ le sous-graphe ayant pour ar\^etes $a_2,\ldots
,a_n$, alors $\P(\cal{G},C)$ est l'amalgame de $G_{s_1}$ et de
$\P(\cal{G},C_1)$ le long de $\varphi_{a_1}$. Alors $u\in
G_{s_1}\cap \P(\cal{G},C_1)$, et donc $u\in G_{a_1}^-\subset
G_{s_1}$ et $u\in G_{a_1}^+\subset G_{\mathrm{e}(a_1)}\subset
\P(\cal{G},C_1)$, ce qui fournit le trajet en $u$, d'origine $s_1$
et d'extr\'emit\'e $\mathrm{e}(a_1)$
$$u=u\overset{a_1}{\longrightarrow} u=u$$
Maintenant, dans $\P(\cal{G},C_1)$, $u\in G_{\mathrm{e}(a_1)}\cap
G_{s_2}$. On proc\`ede au m\^eme raisonnement dans
$\P(\cal{G},C_1)$, et ainsi de suite. Le trajet produit des
trajets successivement
d\'etermin\'es, est le trajet souhait\'e. \hfill$\blacksquare$\\

\noindent On obtient imm\'ediatement le corollaire suivant :

\begin{cor} Soit $(\cal{G},X,T)$
un graphe de groupe d\'ecompos\'e. Soient $s_1,s_2,s_3,s_4$ des
sommets de $X$, et des \el s  $u,v\in \P(\cal{G},X)$, tels que
$u\in G_{s_1}\cap G_{s_2}$, et $v\in G_{s_3}\cap G_{s_4}$.

Alors il existe un trajet de $u$ \`a $v$ d'origine $s_1$ et
d'extr\'emit\'e $s_3$, si et seulement si il existe un trajet de
$u$ \`a $v$, d'origine $s_2$ et d'extr\'emit\'e $s_4$ ; de plus,
ils ont m\^eme label.
\end{cor}

\begin{defn}
\label{C(u)} 
Ainsi $\cal{C}(u,s)$ ne d\'epend pas du choix de $s$, et on le
notera $\cal{C}(u)$. De m\^eme
 on parlera  de
trajet de $u$ \`a $v$, et de circuit en $u$.
\end{defn}

Le r\'esultat qui suit, motive l'introduction des notions de
circuits et de trajets. C'est le r\'esultat fondamental de cette
section.

\begin{prop}[{\bf Th\'eor\`eme fondamental des trajets}]
\label{trajet}
Soit $(\cal{G},X,T)$ un \linebreak graphe de groupe d\'ecompos\'e.
Soient les \el s $u,v$ de $\P(\cal{G},X)$ dans les sous-groupes de
sommet respectifs $G_{s_1},G_{s_2}$.
Alors $u=hvh^{-1}$ dans $\P(\cal{G},X)$, si et seulement si il
existe un trajet de $u$ \`a $v$ ayant pour label $h$.
\end{prop}

\noindent\textbf{D\'emonstration.} Nous avons d\'ej\`a vu
(d\'efinition \ref{traj}), que s'il existe un trajet de $u$ \`a
$v$ ayant pour label $h$ alors $u=hvh^{-1}$ ; montrons la
r\'eciproque.
Soient donc, $u\in G_{s_1}$, et $v\in G_{s_2}$ tels que
$u=hvh^{-1}$ dans $\P(\cal{G},X)$. Soit l'ar\^ete $a\in \cal{A}_X$
; on d\'ecompose $(\cal{G},X)$ le long de $a$ afin d'obtenir le
sous-graphe $(\cal{G},X_1)$.

Si $a$ est non $T$-s\'eparante, $\P(\cal{G},X)$ est l'extension
HNN de $\P(\cal{G},X_1)$ le long de $\varphi_a:
G_a^-\longrightarrow G_a^+$. Les \el s $u$ et $v$ sont dans
$\P(\cal{G},X_1)$, et on peut  donc les supposer cycliquement
r\'eduits de longueur 1. Puisque $u=hvh^{-1}$ dans
$\P(\cal{G},X)$, avec le corollaire \ref{hnnconj2}, soit $h$
conjugue $v$ en $u$ dans $\P(\cal{G},X_1)$, soit $h$ s'\'ecrit
sous forme r\'eduite $h=h_1t_a^{\e_1}\cdots t_a^{\e_p}h_{p+1}$, et
il existe une suite finie  $\g_1,\ldots ,\g_m$ de $G_a^-\cup
G_a^+$,
 avec $m=2p$ pour $p\in \N_\ast$,  v\'erifiant les conditions :
$$u=h_1 \g_1 h_1^{-1} \qquad \g_m=h_{p+1} v h_{p+1}^{-1}\qquad \mathrm{dans\;}
 \P(\cal{G},X_1)$$
$$\g_2=\varphi_a^{\e_1}(\g_1)$$
et pour tout $i=2,\ldots ,p$,
\begin{align*}
\g_{2i-1} &= h_i \g_{2(i-1)} h_i^{-1}\quad \mathrm{dans\; } \P(\cal{G},X_1)\\
\g_{2i} &= \varphi_a^{\e_i}(\g_{2i-1})\qquad \e_i=\pm 1
\end{align*}
Dans ce cas, on symbolise ces conditions par la notation suivante,
que l'on nomme un pr\'e-trajet de $u$ \`a $v$,

$$u\underset{h_1}{\sim} \g_1
\overset{\e_1a}{\longrightarrow} \g_2\underset{h_2}{\sim}
\g_3\overset{\e_2a}{\longrightarrow} \cdots
\overset{\e_pa}{\longrightarrow} \g_m \underset{h_{p+1}}{\sim} v$$
avec $-1\, a$ qui d\'enote $-a$, et $1\, a$ qui d\'enote $a$.
L'\el\ $h_1t_a^{\e_1}\cdots t_a^{\e_p}h_{p+1}$ de $\P(\cal{G},X)$
est appel\'e label du
 pr\'e-trajet. Il faut
remarquer qu'avec cette notation,
$\g\overset{a}{\longrightarrow}\g'$, implique que $\g\in G_a^-$ et
$\g'\in G_a^+$.

Dans le cas o\`u $h$ est de longueur $1$, on consid\`ere le
pr\'e-trajet de $u$ \`a $v$ :
$$u\underset{h}{\sim}v$$

Si $a$ est $T$-s\'eparante, $\P(\cal{G},X)$ est l'amalgame des
facteurs $\P(\cal{G},X_1)$ et $\P(\cal{G},X_2)$ le long de
$\varphi_a:G_a^-\longrightarrow G_a^+$. Avec le corollaire
\ref{amalgconj2}, soit $h$ conjugue $v$ en $u$ dans un des
facteurs, soit $h$ s'\'ecrit sous forme r\'eduite $h=h_1\cdots
h_{m+1}$, et il existe une suite  $\g_1,\ldots ,\g_m$ de
$G_a^-\cup G_a^+$, o\`u deux \el s successifs de $u,\g_1,\ldots
\g_m,v$ sont conjugu\'es dans un facteur, et deux couples
successifs d'\el s conjugu\'es sont dans des facteurs distincts.
De plus,  on a les $m+2$ \'egalit\'es suivantes,
\begin{gather*}
\g_m=h_{p+1}vh_{p+1}^{-1}\qquad u=h_1\g_1h_1^{-1}\\
\forall i=1,\ldots ,m-1,\qquad \g_i=h_i\g_{i-1}h_i^{-1}
\end{gather*}
 chacune dans un facteur $\P(\cal{G},X_1)$ ou
$\P(\cal{G},X_2)$. Puisque $a$ est $T$-s\'eparante, si
$\mathrm{o}(a)$ et dans $X_1$, alors $\mathrm{e}(a)$ est dans
$X_2$ et inversement. Supposons par exemple, que $\mathrm{o}(a)\in
X_1$. On se donne deux couples successifs, $\g_{i-1},\g_i$
conjugu\'es par $h_{i-1}$ dans un facteur, et $\g_i,\g_{i+1}$
conjugu\'es par $h_i$ dans l'autre facteur. Supposons que
$\g_{i-1},\g_i,h_{i-1}$ soient dans $\P(\cal{G},X_1)$.
 On note alors
$$\g_{i-1}\underset{h_{i-1}}{\sim}\g_i
\overset{a}{\longrightarrow} \g_i\underset{h_i}{\sim}\g_{i+1}$$ Il
faut remarquer que l'on a bien $\g_i\in G_a^-=G_a^+$ et
$\g_i=\varphi_a(\g_i)$. Dans l'autre cas, on note
$$\g_{i-1}\underset{h_{i-1}}{\sim}\g_i
\overset{-a}{\longrightarrow} \g_i\underset{h_i}{\sim}\g_{i+1}$$
et l'on a $\g_i\in G_{-a}^-=G_{-a}^+$ et
$\g_i=\varphi_{-a}(\g_i)$. On symbolise ainsi ces conditions par
la notation que l'on baptise encore de pr\'e-trajet de $u$ \`a
$v$,
$$u\underset{h_1}{\sim} \g_1
\overset{\e a}{\longrightarrow} \g_1\underset{h_2}{\sim}
\g_2\overset{-\e a}{\longrightarrow} \cdots \overset{\e(-1)^{m+1}
a}{\longrightarrow} \g_m \underset{h_{m+1}}{\sim} v$$ de label
$h=h_1\cdots h_{m+1}$.

Si $h$ est de longueur 1, on a le pr\'e-trajet trivial :
$$u\underset{h}{\sim}v$$

Revenons \`a la d\'emonstration de la proposition ; on proc\`ede
par induction en d\'ecomposant successivement $X$ le long de ses
ar\^etes. Apr\`es avoir d\'ecompos\'e $X$ le long de $a$, on
obtient l'un des pr\'e-trajets pr\'ec\'edents, que l'on note
$\cal{D}_1$ ; il est soit trivial, soit de la forme suivante :
$$\g_0^+\underset{h_0}{\sim}\g_1^-
\overset{a_1}{\longrightarrow} \cdots
\overset{a_i}{\longrightarrow} \g_i^+\underset{h_i}{\sim}
\g_{i+1}^- \overset{a_{i+1}}{\longrightarrow} \cdots
\overset{a_m}{\longrightarrow} \g_m^+\underset{h_m}{\sim}
\g_{m+1}^-$$ La notation $\g_i^+{\sim} \g_{i+1}^-$ implique la
condition $\g_i^+=h_i\g_{i+1}^-h_i^{-1}$ dans un des facteurs
$\P(\cal{G},X_1)$ ou $\P(\cal{G},X_2)$, et de plus,
$\g_i^+,\g_{i+1}^-$ sont dans des groupes de sommet. On peut
dès-lors leur appliquer le m\^eme proc\'ed\'e dans
$\P(\cal{G},X_1)$ (ou $\P(\cal{G},X_2)$), en d\'ecomposant $X_1$
(ou $X_2$) le long d'une ar\^ete. On obtient comme
pr\'ec\'edemment un pr\'e-trajet $\cal{E}_i$ de $\g_i^+$ \`a
$\g_{i+1}^-$, de label $h_i$. On construit alors un pr\'e-trajet
$\cal{D}_2$ de $u$ \`a $v$ en substituant dans $\cal{D}_1$ pour
tout $i=0,\ldots ,m$, $\g_i^+{\sim} \g_{i+1}^-$ par le
pr\'e-trajet $\cal{E}_i$. Il est important de remarquer que
$\cal{D}_1$ et $\cal{D}_2$ ont m\^eme label $h$. On proc\`ede de
la m\^eme fa\c con lorsque $\cal{D}_1$ est trivial.

En r\'ep\'etant ce proc\'ed\'e,
 on finit par obtenir le
pr\'e-trajet $\cal{C}$, de label $h$,
$$u\underset{k_0}{\sim}{d_1^-}
\overset{\a_1}{\longrightarrow} \cdots
\overset{\a_i}{\longrightarrow} d_i^+
\underset{k_i}{\sim}d_{i+1}^- \overset{\a_{i+1}}{\longrightarrow}
\cdots \overset{\a_q}{\longrightarrow} d_q^+\underset{h_{q}}{\sim}
v$$ v\'erifiant en outre la condition que les \'egalit\'es
$u=h_1d_1^-h_1^{-1}$, $d_i^+=h_{i+1}d_{i+1}^-h_{i+1}^{-1}$ et
$d_q^+=h_{q+1}vh_{q+1}^{-1}$, ont lieu dans des sous-groupes de
sommet. On a ainsi r\'euni toutes les conditions pour avoir un
trajet de $u$ \`a $v$, \`a l'exception d'une : que $(\a_1,\ldots
\a_q)$ soit un chemin de $X$. A priori rien n'assure que ce soit
le cas, et une manipulation simple montre qu'il est en effet
facile d'avoir $\mathrm{e}(\a_i)\not= \mathrm{o}(\a_{i+1})$.
N\'eanmoins, $d_i^+\in G_{\a_i}^+\subset G_{\mathrm{e}(\a_i)}$,
$d_{i+1}^-\in G_{\a_{i-1}}^-\subset G_{\mathrm{o}(\a_{i+1})}$, et
$d_i^+,d_{i+1}^-$ sont dans un m\^eme groupe de sommet $G_s$. Avec
la proposition \ref{III1}, on peut construire un trajet
$\cal{C}_i^+$, de label $1$, de $d_i^+$ \`a $d_i^+$, d'origine
$\mathrm{e}(\a_i)$ et d'extr\'emit\'e $s$. de m\^eme on construit
un trajet $\cal{C}_{i+1}^-$ de $d_{i+1}^-$ \`a
 $d_{i+1}^-$ d'origine $s$ et d'extr\'emit\'e $\mathrm{o}(\a_{i+1})$.
On substitue alors dans $\cal{C}$,
$$\cdots
\overset{\a_i}{\longrightarrow} d_i^+
\underset{k_i}{\sim}d_{i+1}^- \overset{\a_{i+1}}{\longrightarrow}
\cdots$$ par
$$
\cdots \overset{\a_i}{\longrightarrow} \underbrace{d_i^+=d_i^+
\overset{\b_1}{\longrightarrow} \cdots
\overset{\b_r}{\longrightarrow} d_i^+}_{\cal{C}_i^+}
\underset{k_i}{\sim} \underbrace{d_{i+1}^-
\overset{\d_1}{\longrightarrow} \cdots
\overset{\d_s}{\longrightarrow}
d_{i+1}^-=d_{i+1}^-}_{\cal{C}_{i+1}^-}
\overset{\a_{i+1}}{\longrightarrow} \cdots
$$
Cette op\'eration n'a pas modifi\'e le label du pr\'e-trajet, et
cette fois-ci,\linebreak $(\a_i,\b_1,\ldots ,\b_r,\d_1,\ldots
\d_s,\a_{i+1})$ est un chemin. On proc\`ede ainsi pour tout
$i=1,\ldots ,q$, pour peu que ce soit n\'ecessaire. On peut alors
changer la notation $\sim$ par $\circlearrowleft$, car le
pr\'e-trajet obtenu, est un trajet de $u$ \`a $v$ dans
$\P(\cal{G},X)$, ayant pour label $h$, ce qui conclut la
preuve.\hfill$\blacksquare$

On obtient immédiatement :
\begin{cor}
\label{circuit}
Soit $(\cal{G},X,T)$ un graphe de groupe d\'ecompos\'e. Soit $G_s$
un sous-groupe de sommet de $\P(\cal{G},X)$, et $u\in G_s$. Alors
le centralisateur $\cal{Z}(u)$ de $u$ dans $\P(\cal{G},X)$ est
l'ensemble $\cal{C}(u)$ des labels des circuits en $u$.
\end{cor}


\subsection{Th\'eor\`eme de commutativit\'e}

Avec l'\'etude faite dans la section pr\'ec\'edente, nous pouvons
g\'en\'eraliser les th\'eor\`e\-mes \ref{amalgcom} et
\ref{hnncom}, qui caract\'erisent les \el s qui commutent dans un
produit amalgam\'e et dans une extension HNN, au cas du groupe
fondamental d'un graphe de groupe.

\begin{thm}[{\bf Th\'eor\`eme de commutativit\'e dans un graphe de groupe}]
\label{graph_com}
Soit $(\cal{G},X)$ un graphe de groupe muni d'une d\'ecomposition.
Soient $\G=\P (\cal{G},X)$, et $x,y\in \G$ des \'el\'ements qui
commutent. Alors,
\begin{itemize}
\item[(i)] Si $x$ est dans un sous-groupe d'ar\^ete, alors
$y\in \cal{C}(x)$, \emph{i.e.} $y$ est le label d'un circuit en
$x$.
\item[(i')] x ou y est dans le conjugu\'e d'un sous-groupe d'ar\^ete.
\item[(ii)] Si $x$ est dans le conjugu\'e d'un sous-groupe de sommet $G_s$,
et n'est pas dans le conjugu\'e d'un sous-groupe d'ar\^ete, alors
$y$ est dans le m\^eme conjugu\'e de $G_s$.
\item[(iii)] Si ni $x$ ni $y$ n'est dans le conjugu\'e d'un sous-groupe
de sommet, alors $x=ghg^{-1}\, W^j$, $y=gh'g^{-1}\, W^k$, o\`u
$g,W\in \G$, $h,h'$ sont dans un sous-groupe d'ar\^ete $G_a^-$, et
$ghg^{-1},gh'g^{-1},W$ commutent deux \`a deux.
\end{itemize}
\end{thm}
\noindent \textbf{D\'emonstration.} \textsl{Cas} (i). C'est le
corollaire \ref{circuit}. \hfill$\square$

\noindent \textsl{Cas} (i'). Il n'y a rien a montrer dans ce
cas.\hfill$\square$

\noindent \textsl{Cas} (ii). Si $x$ est dans un groupe de sommet,
et n'est pas dans le conjugu\'e d'un groupe d'ar\^ete, tout
circuit en $x$ est trivial, et donc le corollaire \ref{circuit}
permet de conclure.
Si $x$ est dans  le conjugu\'e d'un groupe de sommet $g G_s
g^{-1}$,  alors $g^{-1} x g$ est dans $G_s$, $g^{-1} x g$ et
$g^{-1} y g$ commutent, et le m\^eme raisonnement s'applique.
\hfill $\square$

Pour achever la d\'emonstration, on proc\`ede par induction sur le
nombre d'ar\^etes de $X$. Si $X$ n'a qu'une ar\^ete, les
th\'eor\`emes \ref{amalgcom} et \ref{hnncom} permettent de
conclure. Il est essentiel de remarquer que pour tout $g\in \G$,
si $x,y\in \G$ commutent et v\'erifient une des conclusions (i),
(i'), (ii) ou (iii), alors, $gxg^{-1}$ et $gyg^{-1}$ commutent et
v\'erifient (i), (i'), (ii) ou (iii). Plus pr\'ecis\'ement,
lorsque $x$ et $y$ v\'erifient (i), (i'), (ii), (iii), $gxg^{-1}$
et $gyg^{-1}$ v\'erifient respectivement (i'), (i'), (ii), (iii).

On consid\`ere une ar\^ete $\a$ de $X$. On d\'ecompose $X$ le long
de $\a$. Si $\a$ est non $T$-s\'eparante, on note $Y$ le graphe
obtenu. Alors $\G=\P (\cal{G},X)$ est une extension HNN de $\G'=\P
(\cal{G},Y)$. Avec le th\'eor\`eme \ref{hnncom}, soit on a la
conclusion (i) (i') ou (iii), et l'on peut conclure, soit il
existe $g\in \G$, tel que $x$ et $y$ sont dans $g^{-1} \G' g$. On
pose alors $x'=g x g^{-1}$ et $y'=g y g^{-1}$. Les \'el\'ements
$x'$ et $y'$ sont dans $\G'$ et l'on peut appliquer l'induction.

Si $\a$ est $T$-s\'eparante, on note $Y_1$ et $Y_2$ les deux
composantes connexes du graphe obtenu en d\'ecomposant $Y$ le long
de $\a$. Alors $\G$ est un produit amalgam\'e de $\G_1=\P
(\cal{G},Y_1)$ et $\G_2=\P (\cal{G},Y_2)$. Ainsi, avec le
th\'eor\`eme \ref{amalgcom}, soit $x$ et $y$ v\'erifient les
conditions (i), (i') ou (iii), soit il existe $g\in \G$, tel que
$x$ et $y$ soient dans $g^{-1} \G_1 g$ (respectivement $g^{-1}
\G_2 g$). On pose alors $x'=gxg^{-1}$ et $y'=gyg^{-1}$, $x',y'\in
\G_1$ (respectivement $x',y'\in \G_2$), et l'on peut appliquer
l'induction.\hfill $\blacksquare$

\subsection{Th\'eor\`eme de conjugaison}

Nous nous int\'eressons maintenant au cas de deux \'el\'ements
conjugu\'es dans un graphe de groupe. Ce résultat n'est cependant
qu'une première approximation : il peut être raffiné en un
résultat plus fort, analogue du cas des amalgames (comme montré
dans un cadre plus restreint dans \cite{cpgog3m}) ; c'est l'objet
d'un travail de l'auteur à ce jour en préparation.

Nous avons d'abord besoin d'introduire une certaine proc\'edure de
re-\'ecriture de mots.

\begin{defn}
 Soit $(\cal{G},X,T)$  un
graphe de groupe d\'ecompos\'e. Un \textsl{ordre de
d\'ecom\-position}  est un ordre (au sens large, \emph{i.e.} une
relation r\'eflexive, antisym\'e\-trique et transitive) total
$\preceq$ sur $\cal{A}^+_X$, tel que
$$\forall\, (\a,\b) \in (\cal{A}^+_X-\cal{A}^+_T) \times \cal{A}_T^+\, ,
\quad \a \preceq \b$$
\end{defn}

Intuitivement, un ordre de d\'ecomposition, n'est rien d'autre que
le choix d'un ordre total sur $\cal{A}^+_X$, pour lequel
 les ar\^etes non $T$-s\'eparantes
pr\'ec\`edent les ar\^etes $T$-s\'eparantes. Puisque $\cal{A}_X^+$
est un ensemble fini, tout sous-ensemble $E$ non vide admet un
minimum pour $\preceq$, que nous noterons $\mathrm{min}_\preceq
(E)$ ou plus simplement $\mathrm{min} (E)$. Ainsi on peut parler
du $1^{er},2^{\grave{e}me},\ldots ,p^{i\grave{e}me}$ \el\ de
$\cal{A}_X^+$.

Nous allons d\'ecomposer le graphe de groupe le long de toutes ses
ar\^etes. L'ordre de d\'ecomposition  d\'ecrira l'ordre dans
lequel nous effectuerons cette d\'ecomposition. Il sera purement
arbitraire, nous avons juste, par commodit\'e, souhait\'e
d\'ecomposer le $\P$ d'abord en extension HNN, puis en amalgame.

\begin{defn}
 Un sous-graphe de groupe
$(\cal{G},Y)$ de $(\cal{G},X)$, sera dit
$\mathbf{\preceq}$\textsl{-occurent}, si soit
$(\cal{G},Y)=(\cal{G},X)$,  soit il existe un entier $n \geq 1$,
tel que  le graphe $Y$ soit une composante connexe du graphe
obtenu en d\'ecomposant $X$ le long des $n$ premiers \el s de
$\cal{A}_X^+$. C'est \`a dire que c'est un sous-graphe de groupe
apparaissant lors de la d\'ecomposition de $(\cal{G},X)$ impos\'ee
par l'ordre de d\'ecomposition $\preceq$.
\end{defn}

\noindent\textbf{Proc\'edure de r\'eductions cycliques
successives : }\\
Soit un graphe de groupe $(\cal{G},X,T)$, que l'on munit d'un
ordre de d\'ecomposition $\preceq$. Donn\'ee une pr\'esentation
des sous-groupes de sommets de $(\cal{G},X)$, on dispose d'une
pr\'esentation canonique de $\G=\P(\cal{G},X)$, et d'une famille
g\'en\'eratrice not\'ee $\cal{G}en(X)$ (cf. th\'eor\`eme
\ref{graphpres}, remarque 5), que nous fixons dans la suite.

Consid\'erons un \el\ $u\in \G$, donn\'e par un mot sur la famille
g\'en\'eratrice $S=\cal{G}en(X)$ de $\G$.
 Notons $\a_1=\mathrm{min}(\cal{A}_X^+)$ ;
on d\'ecompose le graphe de groupe $(\cal{G},X)$ le long de
l'ar\^ete $\a_1$.

Si $\a_1$ est non $T$-s\'eparante, soit $X_1$ le graphe obtenu en
d\'ecomposant $X$ le long de $\a_1$. Alors $\G$ est l'extension
HNN de $\G_1=\P(\cal{G},X_1)$ le long de
$\vf_{\a_1}:G_{\a_1}^-\longrightarrow G_{\a_1}^+$. Remarquons que
$S=S_1\cup\{ t_{\a_1}\}$, o\`u $S_1=\cal{G}en(X_1)$ est une
famille g\'en\'eratrice de $\G_1$. Au sens de la d\'ecomposition
HNN de $\G$, on consid\`ere un mot $u_1$ sur $S$, cycliquement
r\'eduit repr\'esentant un conjugu\'e
de $u$ dans $\G$.\\
-- Si $|u_1|>1$, alors $u_1$ est un repr\'esentant satisfaisant
de la classe de conjugaison de $u$, et la proc\'edure s'arr\^ete. \\
-- Si $|u_1|=1$, alors $u_1$ est un mot sur $S_1$, et repr\'esente
un \el\ de $\G_1$. On consid\`ere l'ordre de d\'ecomposition
induit par $\preceq$ sur $\cal{A}_{X_1}^+$, et on r\'ep\`ete la
proc\'edure \`a $u_1$ dans $(\cal{G},X_1,T_1)$, (en posant
$T_1=T\cap X_1$).

Si $\a_1$ est $T$-s\'eparante, on consid\`ere les graphes
$X_1,X_2$, obtenus en d\'ecomposant $X$ le long de $\a_1$, et
$\G_1=\P(\cal{G},X_1), \G_2=\P(\cal{G},X_2)$. Alors $\G$ est
l'amalgame des groupes $\G_1$ et $\G_2$ le long de $\vf_{\a_1} :
G_{\a_1}^-\longrightarrow G_{\a_1}^+$. Au sens de cette
d\'ecomposition de $\G$, on peut consid\'erer un mot $u_1$ sur $S$
cycliquement r\'eduit, conjugu\'e \`a $u$ dans $\G$. \\
-- Si $|u_1|>1$, alors $u_1$ est un repr\'esentant de la classe de
conjugaison
de $u$ satisfaisant et la proc\'edure s'arr\^ete.\\
-- Si $|u_1|=1$, alors $u_1$ est un mot sur $\cal{G}en(X_1)$ ou
$\cal{G}en(X_2)$, et respectivement $u_1 \in \G_1$ ou $u_1\in
\G_2$. Sans perte de g\'en\'eralit\'e, supposons que $u_1\in
\G_1$. On consid\`ere l'ordre de d\'ecomposition induit par
$\preceq$ sur $\cal{A}_{X_1}^+$, et on r\'eapplique la m\^eme
proc\'edure \`a $u_1$ dans $(\cal{G},X_1,T_1)$.\smallskip\\
\indent En r\'ep\'etant ce proc\'ed\'e, on finit par trouver un
mot $u_1$, repr\'esentant de la classe de conjugaison de $u$ dans
$\G$, avec un entier $p\geq 0$, et un graphe $X_p$
$\preceq$-occurent, composante connexe du graphe obtenu en
d\'ecomposant $X$ le long des $p$ premiers \el s $\a_1,\a_2,\ldots
,\a_p$ de $\cal{A}_X^+$ (en posant $X_0=X$), tel que
$u_1\in \P(\cal{G},X_p)$, et soit :\\
-- Le graphe $X_p$ est r\'eduit \`a un sommet $s$, $u_1$
est un mot sur $\cal{G}en(s)$, et $u_1\in G_s$.\\
-- Le graphe $X_p$ n'est pas r\'eduit \`a un sommet. Lorsque l'on
d\'ecompose $X_p$ le long de l'ar\^ete
$\a_{p+1}=\mathrm{min}(\cal{A}_{X_p}^+)$, $\P(\cal{G},X_p)$ se
d\'ecompose en une extension HNN ou un amalgame, et au sens de
cette d\'ecomposition, $u_1$ est un mot cycliquement r\'eduit (sur
$\cal{G}en(X_p)$), de longueur $|u_1|>1$.

\begin{thm}[\textbf{Th\'eor\`eme de conjugaison dans un graphe de groupe}]
\label{graph_conj}
 \hfill \linebreak
Soit
 $(\cal{G},X,T)$ un graphe de groupe d\'ecompos\'e muni
d'un ordre de d\'ecomposition $\preceq$. \linebreak[4] Soient
$u,v\in \P(\cal{G},X)$ des \el s conjugu\'es, et les mots
$u_1,v_1$ obtenus en appliquant la proc\'edure de r\'eductions
cycliques successives \`a des mots sur $\cal{G}en(X)$
repr\'esentant $u$ et $v$.
Alors soit :\medskip\\
\emph{(i)} Il existe deux sommets $s,s'$ de $X$, tels que $u_1$ et
$v_1$ soient des  mots respectivement
 sur $\cal{G}en(s)$ et $\cal{G}en(s')$.
En particulier, $u_1$ et $v_1$ repr\'esentent dans $\P(\cal{G},X)$
des \el s des sous-groupes de sommet $G_s,G_{s'}$, et il existe un
trajet de $u_1$
\`a $v_1$ dans $(\cal{G},X,T)$.\medskip\\
\emph{(ii)} $u_1$ et $v_1$ sont des mots sur $\cal{G}en(X_p)$ o\`u
$X_p$ est un sous-graphe $\preceq$-occurent de $X$, et
repr\'esentent des \el s conjugu\'es dans $\G_p=\P(\cal{G},X_p)$.
En d\'ecomposant $X_p$ le long de l'ar\^ete
$\a_{p+1}=\mathrm{min}(\cal{A}_{X_p}^+)$, $\G_p$ se d\'ecompose,
et au sens de cette d\'ecomposition, $u_1,v_1$ sont de m\^eme
longueur $|u_1|=|v_1|>1$ \\
(les notations sont celles du
paragraphe pr\'ec\'edent).
\end{thm}

\noindent\textbf{D\'emonstration.} Consid\'erons donc sous ces
hypoth\`eses, deux \el s $u,v$ de $\G=\P(\cal{G},X)$ conjugu\'es,
donn\'es par des mots sur la famille g\'en\'eratrice
$\cal{G}en(X)$. On leur applique la proc\'edure de r\'eductions
cycliques successives, pour obtenir des
 mots $u_1,v_1$, repr\'esentants de leur classe de conjugaison
dans $\G$. Observons de plus pr\`es, ce qui peut se passer \`a
chaque \'etape de la d\'ecomposition de $X$. On d\'ecompose $X$ le
long d'une ar\^ete $\a$. Le groupe $\G$ se d\'ecompose soit en une
extension HNN de $\G_1$, soit en un amalgame de $\G_1,\G_2$, le
long de $\vf_\a :G_\a^-\longrightarrow G_\a^+$, selon si $\a$ est
ou non $T$-s\'eparante, et $\G_1$, (resp. et $\G_2$) est(sont)
le(s) groupe(s) du(des) graphe(s) obtenu(s). Avec les corollaires
\ref{amalgconj2} et \ref{hnnconj2}, dans cette d\'ecomposition de
$\G$ on obtient des repr\'esentants cycliquement r\'eduits
$u_0,v_0$ des classes de conjugaison de $u$ et $v$, qui ont m\^eme
longueur,
et l'on est dans exactement un des cas suivants :\\
\textbf{cas 1)} $|u_0|=|v_0|>1$.\\
\textbf{cas 2)}  $|u_0|=|v_0|=1$, et $u_0,v_0$ sont dans un
m\^eme facteur $\G_1$ ou $\G_2$, et conjugu\'es dans ce facteur.\\
\textbf{cas 3)} $|u_0|=|v_0|=1$, et $u_0,v_0$ ne sont pas
conjugu\'es dans un m\^eme facteur. Dans ce cas $u_0$ est dans un
facteur, et conjugu\'e dans ce facteur \`a un \el\ $c_u$ de
$G_\a^{\pm}$, $v_0$ est dans un facteur et conjugu\'e dans ce
facteur \`a un \el\ $c_v$ de $G_\a^{\pm}$, et $c_u$ et $c_v$ sont
conjugu\'es
dans $\G$.\medskip\\
\noindent \textbf{Dans le cas 1)}, on pose $u_1=u_0$, $v_1=v_0$,
et la proc\'edure s'arr\^ete. On vérifie la conclusion (ii)
du th\'eor\`eme.\\
\textbf{Dans le cas 2)}, si $u_0,v_0$ sont dans un facteur
$\P(\cal{G},X_1)$, et si $X_1$ est r\'eduit \`a un sommet $s$,
alors on pose $u_1=u_0$, $v_1=v_0$ ; $u_1$ et $v_1$ sont
conjugu\'es dans $G_s$, et donc il existe un trajet (trivial) de
$u_1$ \`a $v_1$, et la proc\'edure s'arr\^ete. On vérifie la
conclusion (i) du th\'eor\`eme. Sinon, on applique
la m\^eme proc\'edure \`a $u_0$ et $v_0$ dans $(\cal{G},X_1)$.\\
\textbf{Dans le cas 3)}, puisque $c_u$ et $c_v$ sont conjugu\'es
dans $\G$, et sont dans les sous-groupes d'ar\^ete $G_\a^-$ ou
$G_\a^+$, avec la proposition \ref{trajet}, il existe un trajet
$\cal{C}$  de $c_u$ \`a $c_v$ dans $(\cal{G},X)$. Si les facteurs
sont des groupes de sommet, on pose $u_1=u_0$ et $v_1=v_0$, et on
obtient imm\'ediatement l'existence du trajet souhait\'e de $u_1$
\`a $v_1$ dans $(\cal{G},X)$. On vérifie la conclusion (i) du
th\'eor\`eme, et la proc\'edure s'arr\^ete. Sinon, on applique le
m\^eme proc\'ed\'e de r\'eduction cyclique successive, dans un
facteur, \`a $u_0$ et $c_u$ d'une part, et $v_0$ et $c_v$ d'autre
part.

Sans perte de g\'en\'eralit\'e, supposons que $c_u\in G_\a^-$, et
que l'ar\^ete $\a$ a pour origine un sommet du graphe $X_1$.
Ainsi, $G_\a^-\subset \G_1=\P(\cal{G},X_1)$, et $u_0,c_u\in \G_1$.
Puisque $c_u$ est dans le sous-groupe d'ar\^ete $G_\a^-$ de
$(\cal{G},X_1)$, dans tout d\'ecomposition de $(\cal{G},X_1)$ le
long d'une ar\^ete, $c_u$ est un mot de longueur 1. Pour
poursuivre la proc\'edure, on d\'ecompose $(\cal{G},X_1)$ le long
de l'ar\^ete $\b=\mathrm{min}(\cal{A}_{X_1}^+)$, ce qui
d\'ecompose le groupe $\G_1$ en amalgame ou en extension HNN. Dans
cette d\'ecomposition de $\G_1$, $c_u$ est de longueur 1, et
 puisque
$u_0$ et $c_u$ sont conjugu\'es dans $\G_1$, en r\'eduisant
cycliquement
 $u_0$, on obtient un mot $u_0'$ de longueur 1, dans
la classe de conjugaison de $u_0$ et de $c_u$ dans $\G_1$. Ainsi,
on se trouve dans les cas 2) ou 3) figurant ci-dessus. Puisque
pour tout sous-graphe de groupe $(\cal{G},Y)$ de $(\cal{G},X)$, un
\el\ d'un sous-groupe d'ar\^ete $G_\g^-$ de $\P(\cal{G},Y)$ est de
longueur 1 dans toute d\'ecomposition de $\P(\cal{G},Y)$ le long
d'une ar\^ete de $Y$, le m\^eme argument montre qu'en r\'ep\'etant
le m\^eme proc\'ed\'e, on ne se trouvera jamais dans le cas 1)
ci-dessus. Aussi, on finira par d\'eterminer un conjugu\'e $u_1$
de $u$ dans $\G$ dans un sous-groupe de sommet $G_s$, conjugu\'e
\`a $c_u\in G_\a^-$ dans $\G$, et donc un trajet $\cal{C}_u$ dans
$(\cal{G},X)$ de $u_1$ \`a $c_u$. En proc\'edant de la m\^eme fa\c
con avec $c_v$ et $v_0$, on trouvera de m\^eme un \el\ $v_1$ dans
un sous-groupe de sommet, conjugu\'e de $v$, et un trajet
$\cal{C}_v$ de $v_1$ \`a $c_v$ dans $(\cal{G},X)$. Alors, le
trajet produit $\cal{C}_u.\cal{C}.\cal{C}_v^{-1}$ est un trajet de
$u_1$ \`a $v_1$ dans $(\cal{G},X)$. On vérifie alors la conclusion
(i) du th\'eor\`eme. \hfill $\blacksquare$


\subsection{Probl\`eme de la conjugaison et double d'un groupe}

Nous montrons dans cette section que le probl\`eme de la
conjugaison dans un groupe $G$ se r\'eduit au probl\`eme de la
conjugaison dans le double de $G$. Commen\c cons par d\'efinir la
notion de double d'un groupe.

\begin{defn}\label{dble} Soient $G$ un groupe, et $H_1,H_2,\ldots ,H_n$ des
sous-groupes de $G$. Consid\'erons un copie isomorphe $G'$ de $G$,
et un isomorphisme $\vf : G\longrightarrow G'$. Notons
$H_1',H_2',\ldots ,H_n'$ les images respectives de $H_1,H_2,\ldots
,H_n$ par $\vf$. Consid\'erons le graphe de groupe ayant deux
sommets $s,s'$, et $n$ ar\^etes $\a_1,\a_2,\ldots \a_n$, ayant
pour origine $s$ et extr\'emit\'e $s'$, avec $G_s=G$, $G_{s'}=G'$,
et pour $i=1,\ldots ,n$, $G_{\a_i}^-=H_i$, $G_{\a_i}^+=H_i'$, et
$\vf_{a_i}$ est la restriction de $\vf$ \`a $H_i$. Le groupe
fondamental de $(\cal{G},X)$ est appel\'e le \textsl{double} de
$G$ le long des sous-groupes $H_1,H_2,\ldots H_n$, et pourra
\^etre abusivement not\'e 2$G$.
\end{defn}

\begin{thm}
\label{double}
Soient $G$ un groupe et $2G$ son double comme défini ci-dessus.
Alors $G$ se plonge naturellement dans $2G$, et si $u,v\in G$,
alors $u$ et $v$ sont conjugu\'es dans $2G$ si et seulement si ils
sont conjugu\'es dans $G$.
\end{thm}

\noindent\textbf{D\'emonstration.} Reprenons les notations de la
d\'efinition \ref{dble}. La premi\`ere assertion provient
clairement de la d\'efinition. Puisque $G$ se plonge dans $2G$, si
$u$ et $v$ sont conjugu\'es dans $G$, alors ils sont conjugu\'es
dans $2G$ ; montrons la r\'eciproque. Supposons que $u$ et $v$
soient conjugu\'es dans $2G$. Puisque $u$ et $v$ sont dans le
sous-groupe de sommet $G_s=G$, avec le th\'eor\`eme
\ref{graph_conj}, il existe un trajet r\'eduit $\cal{C}$ de $u$
\`a $v$. Si $\cal{C}$ est trivial, alors $u$ et $v$ sont
conjugu\'es dans $G$, aussi on peut supposer que $\cal{C}$ est non
trivial. Puisque $u$ et $v$ sont dans le m\^eme groupe de sommet
$G_s$,  puisque $X$ n'a que deux sommets $s,s'$, et que toute
ar\^ete $\a_i$ a pour origine $s$ et pour extr\'emit\'e $s'$, le
chemin sous-jacent \`a $\cal{C}$ est n\'ecessairement de longueur
paire. De plus, si $p$ est la longueur du chemin, il existe une
application de $\{ 1,\ldots ,p\}$ dans $\{1,\cdots ,n\}$ (on note
$\s_i$ l'image de $i$), telle que le chemin sous-jacent \`a
$\cal{C}$ soit :
$$(\a_{\s_1},-\a_{\s_2},\ldots ,\a_{\s_{2i-1}},-\a_{\s_{2i}},
\ldots ,\a_{\s_{p-1}},-\a_{\s_p})$$ Ainsi, n\'ecessairement,
$\cal{C}$ est de la forme $\cal{C}_0.\cal{C}_1$ avec $\cal{C}_0$ :
$$u\underset{k_1}{\circlearrowleft} u_1^-
\overset{\a_{\s_1}}{\longrightarrow} u_1^+
\underset{h'}{\circlearrowleft} u_2^-
\overset{-\a_{\s_2}}{\longrightarrow} u_2^+ =u_2^+$$ o\`u
$u,u_1^-,u_2^+,k_1\in G$, $u_1^+,u_2^-,h'\in G'$. Puisque
$u_1^+=\vf(u_1^-),u_2^-=\vf(u_2^+)$, et $h'=\vf(h)$ pour un
certain $h\in G$, alors $u_1^-=hu_2^+h^{-1}$ dans $G$. Ainsi on a
le trajet trivial $\cal{D}$ de $u$ \`a $u_2^+$ :
$$u\underset{k_1h}{\circlearrowleft} u_2^+$$
dans $G_s=G$. Consid\'erons le trajet $\cal{D}.\cal{C}_1$. Il va
de $u$ \`a $v$, et a pour chemin sous-jacent :
$$(\a_{\s_3},-\a_{\s_4},\ldots ,\a_{\s_{2i-1}},-\a_{\s_{2i}},
\ldots ,\a_{\s_{p-1}},-\a_{\s_p})$$ qui est de longueur $p-2$.
Ainsi en appliquant le m\^eme argument \`a $\cal{D}.\cal{C}_1$,
puis, successivement, \`a tous les trajets de $u$ \`a $v$ obtenus,
on finit par construire un trajet trivial de $u$ \`a $v$. Ainsi
$u$ et $v$ sont conjugu\'es dans $G$.\hfill$\blacksquare$

\subsection{Propri\'et\'es  d'un graphe de groupe sans circuit}

Apr\`es l'\'etude faite tout au long de cette section, nous
établissons que si un graphe de groupe ne contient pas de circuit
au sens défini ci-dessous, alors la structure de racine, le centre
et les centralisateurs de son groupe fondamental sont dans un sens
triviaux.
 Cette \'etude est
poursuivie, dans un cadre plus large, dans le chapitre 6 de ma
thèse (\cite{phd}) afin d'inclure le cas d'un graphe de groupe
associé à la
décomposition JSJ d'un 3-variété Haken fermée.\\

\begin{defn}
Nous aurons besoin dans la pratique de parler de
\textsl{sous-trajet} d'un trajet
 $\cal{T}$. Si $\cal{T}$ est donn\'e par,
$$u\circlearrowleft c_1^-
\overset{a_1}{\longrightarrow} \cdots
\overset{a_i}{\longrightarrow} c_i^+
\underset{h_i}{\circlearrowleft} c_{i+1}^-
\overset{a_{i+1}}{\longrightarrow} \cdots
\overset{a_n}{\longrightarrow} c_n^+ \circlearrowleft v$$ pour
consid\'erer un sous-trajet de $\cal{T}$, on consid\`ere un
sous-chemin $(a_p,\ldots ,a_q)$ de $(a_1,\ldots ,a_n)$. On
restreint alors le trajet \`a ce sous-chemin. On a plusieurs fa\c
cons de proc\'eder, qui sont d\'esign\'ees par exemple par les
notations :
$$u\circlearrowleft c_1^-
\overset{a_1}{\longrightarrow} \cdots \underbrace{ c_{p-1}^+
\underset{h_{p-1}}{\circlearrowleft} c_{p}^-
\overset{a_{p}}{\longrightarrow} \cdots
\overset{a_q}{\longrightarrow} c_q^+
\underset{h_q}{\circlearrowleft} c_{q+1}^- }_{\cal{D}} \cdots
\overset{a_n}{\longrightarrow} c_n^+ \circlearrowleft v$$
$$u\circlearrowleft c_1^-
\overset{a_1}{\longrightarrow} \cdots \underbrace{ c_{p-1}^+
\underset{h_{p-1}}{\circlearrowleft} c_{p}^-
\overset{a_{p}}{\longrightarrow} \cdots
\overset{a_q}{\longrightarrow} c_q^+ }_{\cal{D}}
\underset{h_q}{\circlearrowleft} c_{q+1}^- \cdots
\overset{a_n}{\longrightarrow} c_n^+ \circlearrowleft v$$
$$u\circlearrowleft c_1^-
\overset{a_1}{\longrightarrow} \cdots c_{p-1}^+
\underset{h_{p-1}}{\circlearrowleft} \underbrace{c_{p}^-
\overset{a_{p}}{\longrightarrow} \cdots
\overset{a_q}{\longrightarrow} c_q^+ }_{\cal{D}}
\underset{h_q}{\circlearrowleft} c_{q+1}^- \cdots
\overset{a_n}{\longrightarrow} c_n^+ \circlearrowleft v$$%
 Qui
correspondent respectivement aux  sous-trajets :
$$c_{p-1}^+ \underset{h_{p-1}}{\circlearrowleft} c_{p}^-
\overset{a_{p}}{\longrightarrow} \cdots
\overset{a_q}{\longrightarrow} c_q^+
\underset{h_q}{\circlearrowleft} c_{q+1}^-$$
$$c_{p-1}^+ \underset{h_{p-1}}{\circlearrowleft} c_{p}^-
\overset{a_{p}}{\longrightarrow} \cdots
\overset{a_q}{\longrightarrow} c_q^+ = c_q^+$$
$$c_{p}^-= c_{p}^-
\overset{a_{p}}{\longrightarrow} \cdots
\overset{a_q}{\longrightarrow} c_q^+ = c_q^+$$
\end{defn}

Nous pourrons aussi d\'esigner un sous-trajet par l'utilisation de
pointill\'es, comme
dans la d\'efinition qui suit.\\


\begin{defn} Un trajet est dit \textsl{r\'eduit}, lorsqu'il ne contient
pas de
 sous-trajet de la forme suivante, avec
$h\in G_a^+=G_{-a}^-$ :
$$\cdots c_1\overset{a}{\longrightarrow}
\underbrace{c_2\underset{h}{\circlearrowleft}
c_3}_{\mathrm{dans\;} G_a^+}
\overset{-a}{\longrightarrow}c_4\cdots$$ Si un trajet $\cal{C}$
n'est pas r\'eduit, on peut proc\'eder \`a la substitution dans
$\cal{C}$, consistant \`a remplacer
$$\cdots
\overset{a_0}{\longrightarrow} u\underset{h_u}{\circlearrowleft}
c_1\overset{a}{\longrightarrow} c_2\underset{h}{\circlearrowleft}
c_3 \overset{-a}{\longrightarrow}c_4
\underset{h_v}{\circlearrowleft} v \overset{a_1}{\longrightarrow}
\cdots$$ par
$$\cdots
\overset{a_0}{\longrightarrow}
u\underset{h_u\varphi_{-a}(h)h_v}{\circlearrowleft}
 v
\overset{a_1}{\longrightarrow} \cdots$$ Une telle op\'eration est
appel\'ee une \textsl{r\'eduction} de $\cal{C}$.
\end{defn}

Il est clair que tout trajet peut \^etre transform\'e par une
suite finie de r\'eductions en un trajet r\'eduit de m\^eme label.
Il est moins clair que l'ordre des r\'eductions n'importe pas,
\emph{i.e.} que le trajet r\'eduit obtenu est unique. Il est
\'el\'ementaire,  de v\'erifier que c'est cependant bien le cas.

\begin{defn}
Un graphe de groupe d\'ecompos\'e est dit \textsl{sans
circuit}, si pour tout $u\not=1$, tout circuit r\'eduit en $u$ est
trivial. Remarquons que l'adjectif {\it d\'ecompos\'e} est ici
redondant.
\end{defn}

\noindent\textbf{Exemple : } Consid\'erons un graphe de groupe
$(\cal{G},X)$, dont les groupes de sommets $G_s$ v\'erifient tous
la propri\'et\'e suivante : si $G_1,G_2,\ldots G_n$ sont les
sous-groupes d'ar\^ete de $G_s$, alors si $i,j=1,2,\ldots ,n$, et
$i\not= j$, aucun \el\ non trivial de $G_i$ n'est conjugu\'e \`a
un \el\ de $G_j$, et si deux \el s $c,c'$ de $G_i$, sont
conjugu\'es par un \el\ $h_i$ dans $G_s$, alors $c=c'$ et $h_i\in
G_s$. Alors $(\cal{G},X)$ est sans circuit. C'est le cas par
exemple pour le graphe associ\'e \`a une d\'ecomposition JSJ d'une
vari\'et\'e Haken  dont toutes les pièces sont des vari\'et\'es
hyperboliques
de volume fini non élémentaires.\\

Le groupe fondamental d'un graphe de groupe sans circuit a des
centralisateurs, un centre, et une structure de racine, triviales
dans un certain sens, comme énoncé ci-dessous.

\begin{thm}
Soit $(\cal{G},X,T)$ un graphe de groupe d\'ecompos\'e sans
circuit. Soit $u\not=1$ un \el\ de $\P(\cal{G},X)$ et $\cal{Z}(u)$
le centralisateur de $u$ dans $\P(\cal{G},X)$.
Alors, soit :
\begin{itemize}\label{centralizer}
\item[(i)] Si $u$ est dans un sous-groupe de sommet $G_s$,
$\cal{Z}(u)$ est le centralisateur de $u$ dans $G_s$.
\item[(ii)] $u$ est dans un conjugu\'e d'un sous-groupe de sommet.
\item[(iii)] Si $u$ n'est pas dans le conjugu\'e d'un sous-groupe
de sommet, alors $\cal{Z}(u)$ est cyclique infini.
\end{itemize}
\end{thm}

\noindent\textbf{D\'emonstration.} Cas (i). Avec le th\'eor\`eme
\ref{graph_com}, si $x$ est dans un groupe de sommet, et si $y$
commute avec $x$, alors $y$ est le label d'un circuit trivial, et
donc $y$ est dans $G_s$.\smallskip\\
\noindent Cas (iii). Supposons que $y$ soit dans $\cal{Z}(u)$, et
que $x$ ne soit pas dans le conjugu\'e d'un facteur. Avec le
th\'eor\`eme \ref{graph_com}, $x=ghg^{-1}W^r$ et $y=gh'g^{-1}W^s$,
et $ghg^{-1},\- gh'g^{-1}, W$ commutent deux \`a deux. Ainsi, soit
$h=h'=1$, soit $g^{-1}Wg$ commute avec un \el\ non trivial d'un
sous-groupe d'ar\^ete $G_a^-$. Ainsi $g^{-1}Wg$ est le label d'un
circuit r\'eduit, et donc puisque $(\cal{G},X)$ est sans circuit,
$g^{-1}Wg$ est dans le groupe de sommet $G_s$, contenant $G_a^-$,
avec $s=\mathrm{o}(a)$. Ainsi $x$ ou $y$ est dans $gG_sg^{-1}$.
Mais avec (i), ceci implique que $x$ est dans $gG_sg^{-1}$, ce qui
est contradictoire. Ainsi, $h=h'=1$, et $x$ et $y$ sont dans le
groupe cyclique engendr\'e par $W$. Pour conclure, $W$ n'est pas
dans le conjugu\'e
d'un sous-groupe de sommet, et est donc sans torsion.\hfill$\blacksquare$\\

\begin{thm}
\label{triv center} Soient $(\cal{G},X,T)$ un graphe de groupe
d\'ecomposé minimal, sans circuit,  $\Gamma$ son groupe
fondamental, et $Z(\G)$ le centre de $\G$. Alors soit :
\begin{itemize}
\item[(i)] $Z(\G)=<1>$,
\item[(ii)] $X$ est réduit à un sommet $s$, et donc $Z(\G)=Z(G_s)$,
\item[(iii)] $X$ est réduit à un sommet $s$ et une arête, avec
$G_s=<1>$, et donc $Z(\G)=\G=\Z$.
\end{itemize}
\end{thm}

\noindent {\bf Démonstration.} Traitons tout d'abord le cas où $X$
contient au moins deux sommets. Alors nécessairement $T$ contient
un sous-graphe $T_0$ constitué de deux sommets $s_0,s_1$ et d'une
arête $a$ d'origine $s_0$ et d'extrémité $s_1$ ; notons $\G_0$ le
groupe fondamental du sous-graphe de groupe associé à $T_0$. Avec
le théorème \ref{graph center} le centre de $\G_0$ contient le
centre de $\G$. Montrons par l'absurde que $Z(G_0)$ est trivial :
soit $\a\not=1$ dans $Z(G_0)$ ; nécessairement $\a$ est dans le
sous-groupe d'arête $G_a$ (théorème \ref{amalgcenter}). Par
hypothèse de minimalité il existe $u\in G_{s_1}\setminus G_a$ et
on alors le circuit réduit non trivial en $\a$ :
$$\a=\a\overset{a}{\longrightarrow}
\a\underset{u}{\circlearrowleft} \a\overset{-a}{\longrightarrow}
\a=\a$$ ce qui est contradictoire. Ainsi on a la conclusion (i).

Considérons maintenant le cas où $X$ contient un unique sommet. La
conclusion (ii) étant évidente supposons en outre que $X$ contient
au moins une arête $a$ ; notons $\G_0$ le sous-groupe de $\G$
associé au sous-graphe de groupe obtenu en supprimant l'arête $a$
de $X$ ; $\G$ est l'extension HNN de $\G_0$ le long de $\vf_a$.
Nécessairement $\mathrm{Fix}\ a=<1>$ ; en effet si $x\in
\mathrm{Fix}\ a$ on a le circuit réduit non trivial en $x$ :
$$x=x\overset{a}{\longrightarrow}x=x$$
et donc par hypothèse $x=1$. De plus nécessairement il n'existe
pas $n>0$, $u\in \G_0$ et $x\not=1\in \G_0$, tels que
$\vf^n(x)=uxu^{-1}$ ; en effet à contrario on aurait le circuit
réduit non trivial en $x$ suivant :
$$x=x
\underbrace{\overset{a}{\longrightarrow}\cdots
\overset{a}{\longrightarrow} \cdots
\overset{a}{\longrightarrow}}_{n\ \mathrm{fois}}
uxu^{-1}\underset{u^{-1}}{\circlearrowleft} x$$ En appliquant le
théorème \ref{hnncenter}, soit $Z(\G)$ est trivial, soit
$\G_0=<1>$ et on obtient alors la condition
(iii).\hfill$\blacksquare$\\

Rappelons qu'un groupe $G$ est dit avoir une \textsl{structure de
racines triviale}  (SRT), si pour tout $g\in G$, l'ensemble
$\{x\in G\, ;\, \exists n\in \Z_\ast ,x^n=g\}$ est inclus dans un
sous-groupe cyclique.

\begin{thm}\label{srt} Soit $(\cal{G},X,T)$ un graphe de groupe d\'ecomposé,
sans circuit. Soient $g$ et $x$ des \el s non
triviaux de $\P(\cal{G},X)$.
Si $x$ est une racine de $g$, alors soit $x$ et $g$ sont dans un
m\^eme conjugu\'e d'un sous-groupe de sommet, soit $x$ est dans le
sous-groupe cyclique infini $\cal{Z}(g)$. En particulier
$\P(\cal{G},X)$ est SRT si et seulement si tous ses groupes de
sommet sont SRT.
\end{thm}

\noindent \textbf{D\'emonstration.} il suffit de remarquer que
sous ces hypoth\`eses, $x$ et $g$ commutent, et d'appliquer le
th\'eor\`eme \ref{centralizer}.
\hfill$\blacksquare$\\

\newpage

\section{{\bf Appendice : d\'emonstration du th\'eor\`eme 5.4}}

Nous montrons  d'abord que $\g$ est conjugu\'e \`a un \el\
cycliquement r\'eduit. Nous proc\'edons par l'absurde. Supposons
que $\g$ ne soit  pas conjugu\'e \`a un \el\ cycliquement
r\'eduit. Soit $K$ la classe des conjugu\'es  de $\g$, et soit
$\mu$ un \'el\'ement de $K$ de longueur minimale  dans $K$.
Puisque $\mu$ n'est pas cycliquement r\'eduit, alors
n\'ecessairement $|\mu |>1$.

L'\el\ $\mu$ s'\'ecrit sous forme r\'eduite
$\mu=\mu_{1}t^{\e_1}\cdots \mu_{n}t^{\e_{n}}\mu_{n+1}$, et quitte
\`a conjuguer $\mu$ par
 $\mu_{n+1}$, on peut supposer que $\mu_{n+1}=1$, c'est \`a dire
que $\mu=\mu_{1}t^{\e_{1}}\cdots \mu_{n}t^{\e_{n}}$.\\
Puisque  $\mu$ n'est pas cycliquement r\'eduit, $|\mu|>2$, et
$t^{\e_{n}}\mu_{1} t^{\e_{1}}$
est un pinch, et alors :\\
\begin{eqnarray*}
\mu_{n}t^{\e_{n}}\; \mu \; t^{-\e_{n}}\mu_{n}^{-1} & = &
\mu_{n}t^{e_{n}}. \mu_{1}t^{\e_{1}}\mu_2\cdots
t^{\e_{n-1}}\mu_{n}t^{\e_{n}}.
t^{-\e_{n}}\mu_{n}^{-1}\\
   & = & (\mu_{n}\f ^{-\e_{n}}(\mu_{1})\,\mu_{2}).\;t^{\e_{2}}\mu_2\cdots
\mu_{n-1}t^{\e_{n-1}}\\
\end{eqnarray*}
En posant $\mu_1 ' = \mu_{n}\f ^{-\e_{n}}(\mu_{1})\,\mu_{2}$
\begin{eqnarray*}
   & = & \mu_{1}'t^{\e_{2}}\cdots \mu_{n-1}t^{\e_{n-1}}\\
\end{eqnarray*}
qui est dans $K$, et de longueur strictement inf\'erieure \`a
$|\mu |$,
ce qui est contradictoire.\hfill $\square$\smallskip\\
\textsl{Cas} (i). \ Si $\g$ est conjugu\'e \`a un \el\
$c\in C_{+1}\cup C_{-1}$.\\
$$\g=h\; c\; h^{-1} \qquad \hbox{o\`{u}} \qquad
h=h_{1}t^{\e_{1}}\cdots h_{p}t^{\e_{p}}h_{p+1}$$
et $h$ est r\'eduit. Nous raisonnons par induction sur $|h|$.\\
Si $|h|=1$, $\g$ est conjugu\'e \`a $c$ par l'\el\ $h=h_1$, qui
est dans $A$, et la conclusion est donc v\'erifi\'ee.

Supposons que $|h|>1$.
\begin{eqnarray*}
\g & = & h \; c \; h^{-1} \\
   & = & h_{1}t^{\e_{1}}\cdots h_{p}t^{\e_{p}}h_{p+1}\; c \; h_{p+1}^{-1}
t^{-\e_{p}}\cdots t^{-\e_1}h_{1}^{-1}\\
\end{eqnarray*}
Le membre de droite n'est pas  r\'eduit, et donc contient un
pinch. Puisque $h$ est r\'eduit, $t^{\e_{p}}\, h_{p+1}\ c\
h_{p+1}^{-1} \, t^{-\e_{p}}$ est un pinch, \emph{i.e.} $h_{p+1} \
c\ h_{p+1}^{-1} \in C_{\e_{p}}$\ On pose $c_{2p}=c$, et
$c_{2p-1}=h_{p+1}\ c\ h_{p+1}^{-1}$. Alors :
\begin{eqnarray*}
\g & = & h_{1}t^{\e_{1}}\cdots h_{p}t^{\e_{p}}c_{2p-1}
t^{-\e_{p}}h_{p}^{-1}\cdots t^{-\e_1}h_{1}^{-1}\\
   & = & h_{1}t^{\e_{1}}\cdots h_{p}\f^{-\e_{p}}(c_{2p-1})h_{p}^{-1}\cdots
t^{-\e_1}{h_{1}}^{-1}\\
\end{eqnarray*}
en posant $c_{2p-2}=\f^{-\e_{p}}(c_{2p-1})\in C_{-\e_{p}}$
\begin{eqnarray*}
\g & = & h_{1}t^{\e_{1}}\cdots t^{\e_{p-1}} h_{p} \ c_{2p-2}
\ h_{p}^{-1}t^{-\e_{p-1}}\cdots t^{-\e_1}h_{1}^{-1}\\
   & = & h'\; c_{2p-2} \; {h'}^{-1}\\
\end{eqnarray*}
en posant $h'= h_{1}t^{\e_{1}}\cdots t^{\e_{p-1}} h_{p}$. On a
obtenu la suite $(c_{2p-2},c_{2p-1},c_{2p}=c)$, d'\el s de
$C_{+1}\cup C_{-1}$, o\`u
\begin{gather*}
c_{2p-2}=\f^{-\e_p}(c_{2p-1})\\
c_{2p-1}=h_{p+1}.c_{2p}.h_{p+1}^{-1}
\end{gather*}

Puisque $\g = h'c_{2p-2}{h'}^{-1}$, o\`u $|h'|<|h|$, on peut
proc\'eder \`a l'induction, et ce faisant on construira la suite
$c_0,c_1,\ldots c_{2p}$, d'\el s de $C_{+1}\cup C_{-1}$,
v\'erifiant
les conclusions du th\'eor\`eme. \hfill $\square$\smallskip\\
\textsl{Cas} (ii).\ Si $\g$ est conjugu\'e \`a un \'el\'ement
$\g'\in A$, et n'est pas conjugu\'e \`a un \'el\'ement de
$C_{+1}\cup C_{-1}$.
\begin{eqnarray*}
\g' & = & h\; \g \; h^{-1} \\
    & = & h_{1}\cdots t^{\e_{n}}h_{n+1}\; \g
\; h_{n+1}^{-1}t^{-\e_{n}}\cdots h_{1}^{-1} \\
\end{eqnarray*}
qui est dans $A$, et donc soit $n=0$, soit l'\'ecriture n'est pas
r\'eduite. Or $h$ est r\'eduit, et donc, si $n \not= 0\;$,
$t^{\e_{n}}\ h_{n+1}\g h_{n+1}^{-1}\ t^{-\e_{n}}$ est un pinch, ce
qui est impossible, puisque $\g$ n'est pas conjugu\'e \`a un \el\
de  $C_{+1}\cup C_{-1}$. Donc $n=0$, c'est \`a dire $\g$ est
conjugu\'e \`a $\g'$
dans $A$.\hfill $\square$ \smallskip\\
\textsl{Cas} (iii). \ Si $\g$ est conjugu\'e \`a un \el\
cycliquement r\'eduit $\g'=u_{1} t^{\mu_1} \cdots u_{m} t^{\mu_m}
$
$$\g=h \g' h^{-1} =h_1 t^{\e_1} \cdots t^{\e_n} h_{n+1} \; \g' \; h_{n+1}
^{-1} t^{-\e_n} \cdots t^{-\e_1} h_1 ^{-1}$$ o\`u h est r\'eduit.
Nous proc\'edons par induction sur $|h|$.

Si $|h|=1$.
$$\g = h_1 \ \g'\ h_1 ^{-1} = h_1 \ u_{1} t^{\mu_1} \cdots u_{m} t^{\mu_m} \ h_1 ^{-1}$$
Le membre de droite est r\'eduit (car $\g'$ est r\'eduit), tandis
que le membre de gauche $\g$ est cycliquement r\'eduit. Alors $h_1
^{-1}$, et donc $h_1$, est dans $C_{\mu_m}$, et la conclusion est
v\'erifi\'ee.

Si $|h|>1$.\\
\begin{eqnarray*}
\g & = & h_1 t^{\e_1} \cdots t^{\e_n} h_{n+1} \; \g'\; h_{n+1}
^{-1} t^{-\e_n}
\cdots t^{-\e_1} h_1 ^{-1} \\
   & = & h_1 t^{\e_1} \cdots t^{\e_n} h_{n+1} \;
u_{1} t^{\mu_1} \cdots u_{m} t^{\mu_m} \; h_{n+1} ^{-1} t^{-\e_n}
\cdots t^{-\e_1} h_1 ^{-1} \\
\end{eqnarray*}
et le membre de droite est de longueur sup\'erieure \`a $|\g |$,
et donc contient un pinch. Puisque $h$ et $\g'$ sont r\'eduits, ce
ne peut-\^etre que $t^{\e_n} h_{n+1} u_{1} t^{\mu_1}$ ou
$t^{\mu_m} h_{n+1} ^{-1} t^{-\e_n}$.

Si $t^{\e_n} h_{n+1}\; u_{1}\; t^{\mu_1}$ est un pinch ; alors
$\e_n = -\mu_1$ et $h_{n+1}\, u_{1} = \a \in C_{\e_n}$,
\begin{align*}
\g &=  h_1 t^{\e_1} \cdots t^{\e_n} h_{n+1} u_{1} t^{\mu_1} \cdots
u_{n} t^{\mu_m} h_{n+1} ^{-1} t^{-\e_n}
\cdots t^{-\e_1} h_1 ^{-1} \\
\intertext{en remplacant $h_{n+1}\, u_{1}$ par $\a$, et
$h_{n+1}^{-1}$ par $u_1 \a ^{-1}$,}
 \g  &= h_1 t^{\e_1} \cdots h_n t^{\e_n} \a  t^{\mu_1} \cdots u_{m}
t^{\mu_m} u_1 \a ^{-1} t^{-\e_n}h_n^{-1} \cdots t^{-\e_1} h_1 ^{-1} \\
  &= h_1 t^{\e_1} \cdots t^{\e_{n-1}} h_n \f^{-\e_n}(\a)u_2 \cdots
u_m t^{\mu_m} u_1 t^{-\e_n} \f^{-\e_n}(\a ^{-1}) h_n ^{-1}
t^{-\e_{n-1}}
\cdots t^{-\e_1}h_1 ^{-1} \\
  &=  h_1 t^{\e_1} \cdots t^{\e_n-1} h_n \f^{-\e_n}(\a)u_2 \cdots
u_m t^{\mu_m} u_1 t^{-\e_n} \f^{-\e_n}(\a)^{-1} h_n ^{-1}
t^{-\e_{n-1}}
\cdots t^{-\e_1}h_1 ^{-1} \\
\intertext{puisque $\e_n = - \mu_1$, }
  \g &= h_1 t^{\e_1} \cdots t^{\e_{n-1}} h_n \f^{\mu_1}(\a)u_2 t^{\mu_2}
\cdots u_m t^{\mu_m} u_1 t^{\mu_1} \f^{\mu_1}(\a)^{-1} h_n ^{-1}
t^{-\e_{n-1}}
\cdots t^{-\e_1}h_1 ^{-1} \\
\intertext{en posant $h_n ' = h_n \f^{\mu_1}(\a)$, }
  \g &= h_1 t^{\e_1} \cdots t^{\e_{n-1}} h_n ' u_2 t^{\mu_2}\cdots
u_m t^{\mu_m} u_1 t^{\mu_1} {h_n '} ^{-1} t^{-\e_{n-1}}
\cdots t^{-\e_1}h_1 ^{-1} \\
\intertext{et en posant $h' = h_1 t^{\e_1} \cdots t^{\e_{n-1},}
h_n '$}
  \g &= h'\;(u_2t^{\mu_2} \cdots u_m t^{\mu_m} u_1 t^{\mu_1})\;{h'} ^{-1}
\end{align*}
Et alors $\g$ s'obtient en conjuguant un conjugu\'e cyclique
de $\g'$ par un \el\ r\'eduit, $h'$, avec $|h'|<|h|$. \\

Si $t^{\mu_m} \; h_{n+1}^{-1} \; t^{-\e_n}$ est un pinch, alors
$h_{n+1} \in C_{\mu_m}$, et $\e_n = \mu_m$.

\begin{align*}
\g &=  h_1 t^{\e_1} \cdots t^{\e_n} h_{n+1}\; u_{1} t^{\mu_1}
\cdots u_{m} t^{\mu_m} \; h_{n+1} ^{-1} t^{-\e_n}
\cdots t^{-\e_1} h_1 ^{-1} \\
\intertext{apr\`es r\'eduction,} \g &=  h_1 t^{\e_1} \cdots
t^{\e_n} h_{n+1}\, u_{1} t^{\mu_1} \cdots u_{m}
\f^{-\mu_m}(h_{n+1})^{-1} h_n ^{-1}
\cdots t^{-\e_1} h_1 ^{-1} \\
\g &=  h_1 t^{\e_1} \cdots h_n \f^{-\e_n}(h_{n+1}) t^{\e_n} u_{1}
t^{\mu_1} \cdots u_{m} \f^{-\mu_m}(h_{n+1}) ^{-1} h_n ^{-1}
\cdots t^{-\e_1} h_1 ^{-1} \\
\intertext{en posant $h_n'=h_n \f^{-\mu_m}(h_{n+1})$,} \g &=  h_1
t^{\e_1} \cdots h_n' t^{\e_n} u_{1} t^{\mu_1} \cdots t^{\mu_{m-1}}
u_{m} {h_n'} ^{-1}
\cdots t^{-\e_1} h_1 ^{-1} \\
\intertext{et $h' = h_1 t^{\e_1} \cdots t^{\e_{n-1}} h_n' $,} \g
&= h'\; t^{\mu_m} u_1 \cdots t^{\mu_{m-1}} u_m \; {h'} ^{-1}
\end{align*}

Remarquons que l'on a $|h'|<|h|$. Seulement $t^{\mu_m} u_1 \cdots
t^{\mu_{m-1}} u_m$ n'est pas un conju\-gu\'e cyclique de $\g'$, au
sens o\`u nous l'entendons. On distingue maintenant deux cas,
selon si $n = 1$,
ou $n > 1$.\\

Si $n=1$. Puisque $\g'$ est cycliquement r\'eduit, l'\el\ $h'\;
t^{\mu_m} u_1 \cdots t^{\mu_{m-1}} u_m \; {h'}^{-1}$ est r\'eduit.
Or il est \'egal \`a $\g$ qui est cycliquement r\'eduit, et donc,
$u_m \, {h'}^{-1} = \b$ est un \el\ de $C_{\mu_{m-1}}$, ainsi,
$$\g = h'\; t^{\mu_m} u_1 \cdots t^{\mu_{m-1}} u_m \; {h'}^{-1}
= \b^{-1} u_m \; t^{\mu_m} u_1 \cdots t^{\mu_{m-1}} \; \b$$ Et
alors $\g$ est conjugu\'e \`a un conjugu\'e cyclique de $\g'$
par un \el\ de $C_{\mu_{m-1}}$.\\

Si $n > 1$.
$$\g =  h_1 t^{\e_1} \cdots t^{\e_{n-1}} h_n' t^{\mu_m}
u_{1} t^{\mu_1} \cdots t^{\mu_{m-1}} u_{m} {h_n'}^{-1}
t^{-\e_{n-1}} \cdots t^{-\e_1} h_1 ^{-1}$$ Le membre de droite est
de longueur sup\'erieure \`a $|\g |$, et  donc contient un pinch.
Puisque $h'$, et $u_{1} t^{\mu_1} \cdots t^{\mu_{m-1}}u_{m}$ sont
r\'eduits, ce ne peut \^etre que, $t^{\e_{n-1}} h_n' t^{\mu_m}$,
ou $t^{\mu_{m-1}} u_{m} {h_n'}^{-1} t^{-\e_{n-1}}$. Or,
$t^{\e_{n-1}} h_n' t^{\mu_m}$ ne peut pas \^etre un pinch. En
effet, si c'est le cas, alors on a $h_n'=h_n \f^{-\mu_m}(h_{n+1})
\in C_{-\mu_m}$, et donc, puisque $\f^{\mu_m}(h_{n+1}) \in
C_{-\mu_m}$, $h_n \in C_{-\mu_m}$. En se rappelant que
$\mu_m=\e_n$, on obtient que $t^{\e_{n-1}} h_n t^{\e_n}$ est un
pinch, ce qui  contredit le fait que $h$ soit r\'eduit.

Ainsi $t^{\mu_{m-1}} u_m {h_n '}^{-1} t^{-\e_{n-1}}$ est un pinch.
C'est \`a dire, $\mu_{m-1}=\e_{n-1}$, et $u_m {h_n'}^{-1} = \b \in
C_{\mu_{m-1}}$.
\begin{align*}
\g &=  h_1 t^{\e_1} \cdots t^{\e_{n-1}} h_n' t^{\mu_m} u_{1}
t^{\mu_1} \cdots t^{\mu_{m-1}} \b t^{-\e_{n-1}}
\cdots t^{-\e_1} h_1 ^{-1} \\
\intertext{en remplacant $h_n '$ par $\b^{-1} u_m$,} \g &=  h_1
t^{\e_1} \cdots t^{\e_{n-1}} \b^{-1} u_m t^{\mu_m} u_{1} t^{\mu_1}
\cdots t^{\mu_{m-1}} \b t^{-\e_{n-1}}
\cdots t^{-\e_1} h_1 ^{-1} \\
\g &= h_1 t^{\e_1} \cdots h_{n-1} \f^{-\e_{n-1}}(\b)^{-1}
t^{\e_{n-1}}u_m t^{\mu_m} \cdots t^{\mu_{m-2}} u_{m-1}
\f^{-\mu_{m-1}}(\b) h_{n-1}^{-1} \cdots t^{-\e_1} h_1 ^{-1} \\
\intertext{avec $\e_{n-1} = \mu_{m-1}$,} \g &= h_1 t^{\e_1} \cdots
h_{n-1} \f^{-\mu_{m-1}}(\b)^{-1}
t^{\e_{n-1}}(u_m t^{\mu_m} \cdots t^{\mu_{m-2}} u_{m-1} t^{\mu_{m-1}})\\
&\qquad\qquad\qquad\qquad\qquad\qquad\qquad\qquad\qquad\qquad
t^{-\e_{n-1}} \f^{-\mu_{m-1}}(\b) h_{n-1}^{-1} \cdots t^{-\e_1} h_1 ^{-1} \\
\intertext{et finalement, en posant $h'' = h_1 t^{\e_1} \cdots
h_{n-1} \f^{-\mu_{m-1}}(\b)^{-1}t^{\e_{n-1}}$,} \g &= h'' (u_m
t^{\mu_m} \cdots t^{\mu_{m-2}} u_{m-1} t^{\mu_{m-1}}) {h''}^{-1}
\end{align*}

et donc $\g$ est conjugu\'e, \`a un conjugu\'e cyclique de $\g'$,
par un \el\ $h''$, r\'eduit, et $|h''|<|h|$.

Ainsi, dans tous les cas, $\g$ s'obtient \`a partir d'un
conjugu\'e cyclique de $\g'$, en conjuguant par un \'el\'ement
r\'eduit, $h'$ ou $h''$, de longueur strictement inf\'erieure \`a
$|h|$, ce qui
nous permet d'appliquer l'induction.\hfill $\blacksquare$ \\

\end{document}